\documentclass{amsart}
\usepackage{amsfonts, amsmath, amsthm}
\usepackage{epsfig}
\usepackage{psfrag}

\newtheorem{pro}{Proposition}[section]
\newtheorem{thm}[pro]{Theorem}
\newtheorem{lem}[pro]{Lemma}
\newtheorem{clm}[pro]{Claim}

\newtheorem{cnj}[pro]{Conjecture}
\newtheorem{cor}[pro]{Corollary}
\newtheorem{quest}[pro]{Question}

\theoremstyle{definition}
\newtheorem{dfn}[pro]{Definition}

\theoremstyle{remark}

\title{Critical Heegaard surfaces} 
\date{December 1, 2000 and, in revised form, August 6, 2001.}
\address{Mathematics Department, University of Illinois at Chicago}
\email{bachman@math.uic.edu}
\author{David Bachman}

\begin{document}
\begin{abstract}
In this paper we introduce {\it critical surfaces}, which are described via a 1-complex whose definition is reminiscent of the curve complex. Our main result is that if the minimal genus common stabilization of a pair of strongly irreducible Heegaard splittings of a 3-manifold is not critical, then the manifold contains an incompressible surface. Conversely, we also show that if a non-Haken 3-manifold admits at most one Heegaard splitting of each genus, then it does not contain a critical Heegaard surface. In the final section we discuss how this work leads to a natural metric on the space of strongly irreducible Heegaard splittings, as well as many new and interesting open questions.
\end{abstract}
\maketitle 

\footnotetext[1]{To appear in Transactions of the AMS.}

\section{Introduction.}

It is a standard exercise in 3-manifold topology to show that every manifold admits Heegaard splittings of arbitrarily high genus. Hence, a ``random" Heegaard splitting does not say much about the topology of the manifold in which it sits. To use Heegaard splittings to prove interesting theorems, one needs to make some kind of non-triviality assumption. The most obvious such assumption is that the splitting is minimal genus. However, this assumption alone is apparently very difficult to use. 

In \cite{cg:87}, Casson and Gordon define a new notion of triviality for a Heegaard splitting, called {\it weak reducibility}. A Heegaard splitting which is not weakly reducible, then, is said to be {\it strongly irreducible}. The assumption that a Heegaard splitting is strongly irreducible has proved to be much more useful than the assumption that it is minimal genus. In fact, in \cite{cg:87}, Casson and Gordon show that in a non-Haken 3-manifold, minimal genus Heegaard splittings {\it are} strongly irreducible. 

The moral here seems to be this: since the assumption of minimal genus is difficult to make use of, one should pass to a larger class of Heegaard splittings, which is still restrictive enough that one can prove non-trivial theorems. 

Now we switch gears a little. It is a Theorem of Riedemeister and Singer (see \cite{am:90}) that given two Heegaard splittings, one can always stabilize the higher genus one some number of times to obtain a stabilization of the lower genus one (see the next section for definitions of these terms). However, this immediately implies that any two Heegaard splittings have a common stabilization of arbitrarily high genus. Hence, the assumption that one has a ``random" common stabilization cannot be terribly useful. What is of interest, of course, is the minimal genus common stabilization. As before though, the assumption of minimal genus has turned out to be very difficult to use. 

In this paper, we propose a new class of Heegaard splittings, which we call {\it critical}, and prove that at least in the non-Haken case, this class includes the minimal genus common stabilizations. As one would hope, the assumption that a splitting is critical is more useful than the assumption that it is a minimal genus common stabilization. 

We define the term {\it critical} via a 1-complex associated with any embedded, separating surface in a 3-manifold, which is reminiscent of the curve complex. After the preliminary definitions we present a section of Lemmas which lead up to the following Theorem:

\medskip
\noindent {\bf Theorem 4.6.} {\it Suppose $M$ is an irreducible 3-manifold with no closed incompressible surfaces, and at most one Heegaard splitting (up to isotopy) of each genus. Then $M$ does not contain a critical Heegaard surface.} 
\medskip

The remainder of the paper is concerned with the converse of this Theorem. That is, we answer precisely when a (non-Haken) 3-manifold {\it does} contain a critical Heegaard surface. 

The main technical theorem which starts us off in this direction is:

\medskip
\noindent {\bf Theorem 5.1.} {\it Let $M$ be a 3-manifold with critical surface, $F$, and incompressible surface, $S$. Then there is an incompressible surface, $S'$, homeomorphic to $S$, such that every loop of $F \cap S'$ is essential on both surfaces. Furthermore, if $M$ is irreducible, then there is such an $S'$ which is isotopic to $S$.}
\medskip

As immediate corollaries to this, we obtain:

\medskip
\noindent {\bf Corollary 5.6.} {\it If $M_1 \# M_2$ contains a critical surface, then either $M_1$ or $M_2$ contains a critical surface.}
\medskip

\noindent {\bf Corollary 5.7.} {\it A reducible 3-manifold does not admit a critical Heegaard splitting.}

\medskip
\noindent {\bf Corollary 5.8.} {\it Suppose $M$ is a 3-manifold which admits a critical Heegaard splitting, such that $\partial M \ne \emptyset$. Then $\partial M$ is essential in $M$.}
\medskip

It is this last corollary which we combine with a considerable amount of new machinery to yield:

\medskip
\noindent {\bf Theorem 7.1.$'$} {\it Suppose $F$ and $F'$ are distinct strongly irreducible Heegaard splittings of some closed 3-manifold, $M$. If the minimal genus common stabilization of $F$ and $F'$ is not critical, then $M$ contains an incompressible surface.}
\medskip

We actually prove a slightly stronger version of this Theorem, that holds for manifolds with non-empty boundary. 

Compare Theorem 7.1.$'$ to that of Casson and Gordon \cite{cg:87}: If the minimal genus Heegaard splitting of a 3-manifold, $M$, is not strongly irreducible, then $M$ contains an incompressible surface. In this light, we see that critical Heegaard surfaces are a natural follow-up to strongly irreducible Heegaard surfaces. 

In the last section, we show how these results lead to a natural metric on the space of strongly irreducible Heegaard splittings of a non-Haken 3-manifold. It is the belief of the author that more information about this space would be of great interest. We also discuss how a better understanding of critical surfaces may help answer several open questions, including the {\it stabilization conjecture}. We end with a few new questions and conjectures generated by this work which the reader may find interesting. 

The author would like to thank several people for their input during the preparation of this paper: the referee, for many helpful suggestions regarding earlier versions; Cameron Gordon, for guidance during the author's thesis work, from which this paper grew; Saul Schleimer, for an introduction to the beauty of the curve complex; and Eric Sedgwick, for helpful conversations regarding Lemma \ref{l:reducesphere}.

\section{Basic Definitions.}
In this section, we give some of the standard definitions that will be used throughout the paper. The expert in 3-manifold theory can easily skip this. 

A 2-sphere in a 3-manifold which does not bound a 3-ball on either side is called {\it essential}. If a manifold does not contain an essential 2-sphere, then it is referred to as {\it irreducible}.

A loop on a surface is called {\it essential} if it does not bound a disk in the surface. Given a surface, $F$, in a 3-manifold, $M$, a {\it compressing disk} for $F$ is a disk, $D \subset M$, such that $F \cap D=\partial D$, and such that $\partial D$ is essential on $F$. If we let $D \times I$ denote a thickening of $D$ in $M$, then to {\it compress $F$ along $D$} is to remove $(\partial D) \times I$ from $F$, and replace it with $D \times \partial I$.  

A {\it compression body} is a 3-manifold which can be obtained by starting with some surface, $F$, forming the product, $F \times I$, attaching some number of 2-handles to $F \times \{1\}$, and capping off any remaining 2-sphere boundary components with 3-balls. The boundary component, $F \times \{0\}$, is often referred to as $\partial _+$. The other boundary component is referred to as $\partial _-$. If $\partial _-=\emptyset$, then we say the compression body is a {\it handlebody}.

A surface, $F$, in a 3-manifold, $M$, is a {\it Heegaard splitting of M}, if $F$ separates $M$ into two compression bodies, $W$, and $W'$, such that $F=\partial _+W=\partial _+ W'$. Such a splitting is {\it non-trivial} if neither $W$ nor $W'$ are products. A {\it stabilization} of $F$ is a new Heegaard splitting which is the connect sum of the standard genus 1 Heegaard splitting of $S^3$ and $F$. Another way to define a stabilization is by ``tunneling" a 1-handle out of $W$, and attaching it to $\partial _+ W'$. If this is done in such a way so as to make the definition symmetric in $W$ and $W'$, then one arrives at a stabilization. The Riedemeister-Singer theorem (see \cite{am:90}) states that given any two Heegaard splittings, $F$ and $F'$, there is always a stabilization of $F$ which is isotopic to a stabilization of $F'$.

\section{The isotopy-invariant disk complex.}

For every surface, $F$, we can define a complex, $G(F)$, as follows: for each isotopy class of essential simple closed curve in $F$, there is a vertex of $G(F)$. There is an edge connecting two vertices if and only if there are representatives of the corresponding equivalence classes that are disjoint. $G(F)$ is the 1-skeleton of a complex which is usually referred to as the {\it curve complex} of $F$, and is an object that many mathematicians have made use of. We now generalize the 1-complex, $G(F)$. 

The new 1-complex we define here will be only for embedded, orientable, compact, separating (but not necessarily connected) surfaces in 3-manifolds. Suppose $F$ is such a surface, separating a 3-manifold, $M$, into a ``red" and a ``blue" side. If $D$ and $D'$ are compressing disks for $F$, then we say $D$ is equivalent to $D'$ if there is an isotopy of $M$ taking $F$ to $F$, and $D$ to $D'$ (we do allow $D$ and $D'$ to be on opposite sides of $F$). Note that this equivalence relation is stronger than $D$ and $D'$ being isotopic in $M$ (rel $F$). Indeed, if there is an isotopy which takes $D$ to $D'$, sweeping out a subset $B \subset M$, then there is a isotopy of $M$, taking $D$ to $D'$, which fixes every point outside of a neighborhood of $B$. 

We now define a 1-complex, $\Gamma (F)$. For each equivalence class of compressing disk for $F$, there is a vertex of $\Gamma (F)$. Two (not necessarily distinct) vertices are connected by an edge if there are representatives of the corresponding equivalence classes on opposite sides of $F$, which intersect in at most a point. $\Gamma (F)$ is thus an ``isotopy-invariant disk complex" for $F$. 

For example, if $F$ is the genus 1 Heegaard splitting of $S^3$, then there is an isotopy of $S^3$ which takes $F$ back to itself, but switches the sides of $F$. Such an isotopy takes a compressing disk on one side of $F$ to a compressing disk on the other. Hence, $\Gamma (F)$ has a single vertex. However, there are representatives of the equivalence class that corresponds to this vertex which are on opposite sides of $F$, and intersect in a point. Hence, there is an edge of $\Gamma (F)$ which connects the vertex to itself. 

We can now exploit this terminology to give concise definitions of some of the standard terms in 3-manifold topology, as well as one important new one. 

\begin{dfn}
We say $F$ is {\it incompressible} if $\Gamma (F)=\emptyset$. If an irreducible 3-manifold contains no incompressible surfaces, then we say it is {\it non-Haken}. Otherwise, it is {\it Haken}.
\end{dfn}

\begin{dfn}
$F$ is {\it reducible} if there are compressing disks on opposite sides of $F$ with the same boundary. The union of these disks is called a {\it reducing sphere} for $F$. 
\end{dfn}

It is a standard exercise to show that if a Heegaard splitting of an irreducible 3-manifold is reducible, then there is a lower genus Heegaard splitting. 

\begin{dfn}
$F$ is {\it strongly irreducible} if there are compressing disks on opposite sides of $F$, but $\Gamma (F)$ contains no edges.
\end{dfn}

In \cite{cg:87}, Casson and Gordon show that if $F$ is a Heegaard splitting of a non-Haken 3-manifold which is not strongly irreducible, then $F$ is reducible. Hence, any minimal genus Heegaard splitting of a non-Haken 3-manifold must be strongly irreducible. 

\begin{dfn}
A vertex of $\Gamma (F)$ is {\it isolated} if it is not the endpoint of any edge.
\end{dfn}

\begin{dfn}
If we remove the isolated vertices from $\Gamma (F)$ and are left with a disconnected 1-complex, then we say $F$ is {\it critical}.
\end{dfn}

Equivalently, $F$ is critical if there exist two edges of $\Gamma(F)$ that can not be connected by a path.

\section{Preliminary Lemmas}

Before proceeding further, we introduce some notation. Suppose $D$ and $E$ are compressing disks on opposite sides of $F$, such that $|D \cap E| \le 1$. Then we denote the edge of $\Gamma (F)$ which connects the equivalence class of $D$ to the equivalence class of $E$ as $D-E$. If $D$ and $D'$ are disks in the same equivalence class, then we write $D \sim D'$. Hence, a chain of edges in $\Gamma (F)$ may look something like
\[D_1-E_2 \sim E_3-D_4 \sim D_5 - E_6 -D_7\]
Many of the proofs of this paper follow by producing such chains. Note that we will always denote disks on the red side of $F$ (i.e. ``red disks") with the letter ``$D$" (usually with some subscript), and blue disks with the letter ``$E$".  

We begin with a simple Lemma.

\begin{lem}
\label{l:3site}
Let $F$ be a critical surface, and suppose $D_1 - E_1$ and $D_2 - E_2$ are edges which lie in different components of $\Gamma (F)$. Then for each red disk, $D$, there is an $i=1$ or $2$ such that $|D \cap E| >1$, whenever $E$ is a blue disk equivalent to $E_i$. Similarly, for each blue disk, $E$, there is an $i$ such that $|D \cap E| >1$, whenever $D$ is a red disk equivalent to $D_i$. 
\end{lem}

\begin{proof} 
If the Lemma is false, then we have the following chain of edges: $D_1 - E_1 - D - E_2 - D_2 $, which clearly contradicts the assumption that the edges $D_1 - E_1$ and $D_2 - E_2$ are in different components of $\Gamma (F)$.
\end{proof}

{\bf Note.} Lemma \ref{l:3site} is closely related to the ``3-site property" of Pitts and Rubinstein (see, for example, \cite{rubinstein:96}). In our language, the 3-site property says that if $D_1 - E_1$ and $D_2 - E_2$ are edges in different components of $\Gamma (F)$, then there cannot be a third edge, $D - E$, where $D$ misses $E_1$ and $E_2$, and $E$ misses $D_1$ and $D_2$. It is trivial to show that this is implied by Lemma \ref{l:3site}. 

The next Lemma is a nice self-contained result about the components of a disconnected critical surface. It is presented only to help the reader get used to our definitions, and will not be used in the remainder of the paper. 

\begin{lem}
\label{l:components}
Suppose $F$ is a critical surface such that every compressing disk for every component of $F$ is isotopic (rel $\partial$) to a compressing disk for $F$, and such that each component of $F$ is separating. Then either 
\begin{itemize}
    \item exactly one component of $F$ is critical, and the rest are incompressible, or
    \item exactly two components of $F$ are strongly irreducible, and the rest are incompressible.
\end{itemize}  
\end{lem}

\noindent {\bf Note.} It is possible to define the terms {\it strongly irreducible} and {\it critical} for non-separating surfaces as well, which makes the statement of Lemma \ref{l:components} a little nicer. We do not bother here, since we will not need such generality. 

\begin{proof}
Suppose that $D_1-E_1$ and $D_2-E_2$ are edges in different components of $\Gamma (F)$. By Lemma \ref{l:3site}, $\partial D_2 \cap \partial E_1 \ne \emptyset$. Hence, $D_2$ and $E_1$ are compressing disks for the same component of $F$, which we refer to as $F_1$. Similarly, $D_1$ and $E_2$ must be compressing disks for the same component, $F_2$. Any compressing disk for any other component would violate Lemma \ref{l:3site}. Hence, any component of $F$ other than $F_1$ or $F_2$ is incompressible.

Now suppose that $F_1=F_2$. If $D_1-E_1$ and $D_2-E_2$ are in the same component of $\Gamma (F_1)$, then there is a chain of edges that connects them. But every compressing disk for $F_1$ is a compressing disk for $F$, so the same chain would connect $D_1-E_1$ and $D_2-E_2$ in $\Gamma (F)$, a contradiction. Hence, $F_1$ is critical. 

If $F_1 \ne F_2$, and $F_2$ is not strongly irreducible, then let $D - E$ be an edge of $\Gamma (F_2)$. Then we have the chain: 
\[ D_1 -E_1 - D-E-D_2-E_2\]
Since this is a contradiction, $\Gamma (F_2)$ contains no edges, and $F_2$ is strongly irreducible. A symmetric argument shows that $F_1$ is strongly irreducible.
\end{proof}

For the proof of the next Lemma, we will need to define a partial ordering on embedded surfaces in a 3-manifold. This ordering will play a more crucial role in the next section. 

\begin{dfn}
\label{d:order}
For any surface, $F$, let $c(F)=\sum \limits _n (2-\chi(F^n))^2$, where $\{F^n\}$ are the components of $F$. If $F_1$ and $F_2$ denote compact, embedded surfaces in a 3-manifold, $M$, 
then we say $F_1 < F_2$ if $c(F_1) < c(F_2)$.
\end{dfn}

Note that this ordering is defined so that if $F_1$ is obtained from $F_2$ by a compression, then $F_1 < F_2$. In fact, we could have used any complexity we wished (and there are several) which induced this partial ordering. The one used here was introduced to the author by Peter Shalen.

\begin{lem}
\label{l:reducesphere}
Suppose $F$ is a Heegaard splitting of an irreducible 3-manifold, $M$, which does not contain any closed incompressible surfaces. If $D_0$ and $E_0$ are red and blue disks such that $D_0 \cap E_0=\emptyset$, then there is a reducing sphere for $F$ which corresponds to an edge of $\Gamma (F)$ in the same component as $D_0-E_0$.
\end{lem}

\begin{proof}
The proof is a restatement of Casson and Gordon's proof from \cite{cg:87}, which says there exists some reducing sphere for $F$. By being careful, we can guarantee that this reducing sphere corresponds to an edge of $\Gamma (F)$ in the same component as $D_0-E_0$. 

Choose non-empty collections of red and blue disks, $\bf D$ and $\bf E$, subject to the following constraints:
\begin{enumerate}
    \item For all $D \in \bf D$ and $E \in \bf E$, $D \cap E=\emptyset$. 
    \item For all $D \in \bf D$ and $E \in \bf E$, $D - E$ is in the same component of $\Gamma (F)$ as $D_0-E_0$.
    \item The surface, $F_0$, obtained from $F$ by compressing along all disks of $\bf D$ and $\bf E$ is minimal (in the sense of Definition \ref{d:order}), with respect to the above two constraints. 
\end{enumerate}

Note that the existence of $D_0$ and $E_0$ is what guarantees that we may find such collections, $\bf D$ and $\bf E$, which are non-empty. 

We claim that each component of the surface, $F_0$, is a 2-sphere. First, let $F_+$ be the surface obtained from $F$ by compressing along all the disks of $\bf D$, and let $F_-$ be the surface obtained by compressing along all disks of $\bf E$. Let $M_+$ be the closure of the component of $M-F_0$ which contains $F_+$, and let $M_-$ be the closure of the other component. Notice that $F_+$ and $F_-$ are Heegaard splittings of $M_+$ and $M_-$. 

If some component of $F_0$ is not a 2-sphere, then it is compressible into $M_+$ or $M_-$ (since $M$ contains no closed incompressible surfaces). Suppose the former. Let $W$ denote the closure of the component of $M_+-F_+$ which contains $F_0$, and let $W'$ denote the other component. By the Lemma of Haken (see, for example, \cite{cg:87}), there is a complete collection, $\bf E'$, of compressing disks for $F_+$ in $W$, and a compressing disk, $\hat D$, for $F_0$ in $M_+$, such that $\hat D \cap F_+$ is a single loop, and $\hat D \cap E'=\emptyset$ for all $E' \in \bf E'$. Notice that $\hat D$ contains a compressing disk, $\hat D'$, for $F_+$ in $W'$. But since $F_+$ is obtained from $F$ by compression along $\bf D$, we can consider $\hat D'$ and $\bf E'$ to be compressing disks for $F$, which are disjoint from all $D \in \bf D$. 

Now, we replace the original collections, $\bf D$ and $\bf E$, with the collections ${\bf D'}= {\bf D} \cup \hat D'$ and $\bf E'$. By construction, for any $D' \in \bf D'$ and $E' \in \bf E'$, $D' \cap E'=\emptyset$. For any $E' \in \bf E'$ and $D \in \bf D$ (and hence in $\bf D'$), $D-E'$ is an edge in $\Gamma (F)$ which is in the same component as $D-E$, for any $E \in \bf E$, and hence in the same component as $D_0-E_0$. Furthermore, compressing $F$ along all disks of $\bf D'$ and $\bf E'$ yields a surface that can be obtained from $F_0$ by compressing along $\hat D'$. But this surface is then smaller than $F_0$, contradicting our minimality assumption. We conclude that $F_0$ is a collection of 2-spheres. 

The remainder of the proof follows that of Rubinstein from \cite{rubinstein:96}. Notice that $F_0$ consists of a subsurface, $P$, of $F$, together with two copies of each disk in $\bf D$ and $\bf E$. These disks are colored red and blue, so we can picture $F_0$ as a collections of spheres with a bunch of red and blue subdisks. Some sphere in this collection, $F_0'$, must have both red and blue subdisks, because otherwise $F$ would not be connected. Let $\gamma$ be a loop on $F_0'$ which separates the red disks from the blue disks. Note that $\gamma \subset P \subset F$. Let $R$ be the component of $F'_0-\gamma$ which contains the red disks, and let $B$ be the other component. Then we can push $P \cap R$ slightly into the red side, and $P \cap B$ into the blue side. This turns $F'_0$ into a sphere that intersects $F$ in the single essential simple closed curve, $\gamma$. Furthermore, $R \cap E=\emptyset$ for all $E \in \bf E$, so the edge $R-B$ is in the same component of $\Gamma (F)$ as $D_0-E_0$. 
\end{proof}

\begin{lem}
\label{l:stab}
Suppose $F$ is a Heegaard splitting of an irreducible 3-manifold, $M$, which does not contain any closed incompressible surfaces. For each edge, $D-E$ of $\Gamma (F)$, there is some edge, $D'-E'$, in the same component as $D-E$, such that $|D' \cap E'|=1$.
\end{lem}

\begin{proof}
Since $D-E$ is an edge of $\Gamma (F)$, either $|D \cap E|=1$, or $D \cap E=\emptyset$. In the former case there is nothing to prove. In the latter case, we may apply Lemma \ref{l:reducesphere} to obtain an edge, $D''-E''$, in the same component as $D-E$, such that $D'' \cup E''$ is a reducing sphere for $F$. By the irreducibility of $M$, $D'' \cup E''$ bounds a ball, $B$. By \cite{waldhausen:68}, we know that $F \cap B$ is standard. So inside $B$ we can find a pair of disks on opposite sides, $D'$ and $E'$, such that $|D' \cap E'|=1$. But then $D' \cap E''=\emptyset$, so the chain $D''-E''-D'-E'$ implies $D'-E'$ is in the same component of $\Gamma (F)$ as $D-E$. 
\end{proof}

\begin{thm}
\label{t:nocrit}
Suppose $M$ is an irreducible 3-manifold with no closed incompressible surfaces, and at most one Heegaard splitting (up to isotopy) of each genus. Then $M$ does not contain a critical Heegaard surface. 
\end{thm}

Examples of such 3-manifolds include $S^3$ and $B^3$ \cite{waldhausen:68}, $L(p,q)$ \cite{bo:83}, \cite{bonahon:83}, and handlebodies \cite{waldhausen:68}.

\begin{proof}
Let $F$ be a Heegaard splitting of $M$, and let $D-E$ and $D'-E'$ denote edges of $\Gamma (F)$. Our goal is to show that these two edges lie in the same component. By Lemma \ref{l:stab}, we may assume that $|D \cap E|=|D' \cap E'|=1$. As $F$ cannot be a torus, the boundary of a neighborhood of $D \cup E$ is a sphere which intersects $F$ in a single essential simple closed curve. This curve divides the sphere into a red disk, $R$, and a blue disk, $B$, such that $\partial R=\partial B$. Similarly, $D' \cup E'$ gives rise to disks, $R'$ and $B'$, such that $\partial R' =\partial B'$. 

Now, $\partial R$ separates $F$ into two components, at least one of which, $T$, is a punctured torus. Furthermore, $(F-T) \cup R$ is a Heegaard surface of lower genus. Similarly, $\partial R'$ separates $F$ into two components, at least one of which, $T'$, is a punctured torus, and $(F-T') \cup R'$ is a Heegaard surface of lower genus. 

By assumption, $(F-T) \cup R$ is isotopic to $(F-T') \cup R'$. This isotopy can be realized as an isotopy of $M$, which takes $(F-T) \cup R$ to $(F-T') \cup R'$. Furthermore, since $B$ and $B'$ are isotopic to $R$ and $R'$, we may assume that the isotopy of $M$ takes $B$ to $B'$, as long as $B$ and $B'$ do not end up on opposite sides. Since $T$ and $T'$ lie inside the balls bounded by $R \cup B$ and $R' \cup B'$ (and in fact, $T \cup R$ and $T' \cup R'$ are Heegaard splittings for these balls), we can also arrange the isotopy so that it takes $T$ to $T'$. What we have now produced is an isotopy of $M$, taking $F$ to $F$, and $R$ to $R'$. This implies the following chain:
\[E-D-B-R \sim R'-B'-D'-E'\]

If, on the other hand, the isotopy of $M$ which takes $(F-T) \cup R$ to $(F-T') \cup R'$ ``switches sides", so that $B$ and $B'$ end up on opposite sides, then there must be an isotopy which takes $(F-T) \cup R$ to $(F-T') \cup B'$, and $(F-T) \cup B$ to $(F-T') \cup R'$. Hence, we conclude that there is an isotopy of $M$ taking $F$ to $F$, and $R$ to $B'$. This implies the chain:
\[E-D-B-R \sim B'-R'-D'-E'\]
\end{proof}

\section{Incompressible Surfaces}
\label{s:incomp}

Our goal in this section is to examine the interplay between critical surfaces and incompressible surfaces. Our main result is the following:

\begin{thm}
\label{t:essential_intersection}
Let $M$ be a 3-manifold with critical surface, $F$, and incompressible surface, $S$. Then there is an incompressible surface, $S'$, homeomorphic to $S$, such that every loop of $F \cap S'$ is essential on both surfaces. Furthermore, if $M$ is irreducible, then there is such an $S'$ which is isotopic to $S$.
\end{thm}

\begin{proof}
The proof is in several stages. First, we construct a map, $\Phi$, from $S \times D^2$ into $M$. We then use $\Phi$ to break up $D^2$ into regions, and label them in such a way so that if any region remains unlabelled, then the conclusion of the theorem follows. Finally, we construct a map from $D^2$ to a labelled 2-complex, $\Pi$, which has non-trivial first homology, and show that if there is no unlabelled region, then the induced map on homology is nontrivial, a contradiction. This general strategy is somewhat similar to that used in \cite{rs:96}, although the details have very little in common. 

\medskip
\noindent \underline{Stage 1:} {\it Constructing the map, $\Phi :S \times D^2 \rightarrow M$}.
\medskip

Let $D_0-E_0$ and $D_1-E_1$ denote edges in different components of $\Gamma (F)$, such that $D_0 \cap E_0=D_1 \cap E_1=\emptyset$. Such edges are guaranteed to exist in every component. Indeed, if $D-E$ is an edge where $|\partial D \cap \partial E|=1$, then the boundary of a neighborhood of $D \cup E$ is a sphere, which intersects $F$ in a single loop. This loop cuts the sphere into disks, $D'$ and $E'$, which can be isotoped to be disjoint. Furthermore, we have the chain $D'-E'-D-E$ in $\Gamma (F)$, which shows that $D'-E'$ is in the same component as $D-E$. 

We now produce a sequence of compressing disks for $F$, $\{D_{\frac {i}{n}}\}_{i=1,...,n-1}$, such that $D_{\frac {i}{n}} \cap D_{\frac {i+1}{n}} =\emptyset$, for all $i$ between 0 and $n-1$. Begin by isotoping $D_0$ and $D_1$ so that $|D_0 \cap D_1|$ is minimal. If $D_0 \cap D_1=\emptyset$, then we are done. If not, then let $\gamma$ be some arc of intersection which is outermost on $D_1$. So, $\gamma$ cuts off a subdisk, $D_1'$, of $D_1$, whose interior is disjoint from $D_0$. $\gamma$ also cuts $D_0$ into two subdisks. Choose one to be $D_0'$. Since $|D_0 \cap D_1|$ is assumed to be minimal, the disk $D^1=D_1' \cup D_0'$ must be a compressing disk for $F$. Furthermore, $D^1$ can be isotoped so that is misses $D_0$, and so that $|D^1 \cap D_1|$ is strictly less than $|D_0 \cap D_1|$. Continue now in the same fashion, using $D^1$ and $D_1$ to construct $D^2$, where $D^1 \cap D^2 =\emptyset$, and $|D^2 \cap D_1|<|D^1 \cap D_1|$. Eventually, we come to a disk, $D^{n-1}$, where $D^{n-1} \cap D_1=\emptyset$. Now, for all $i$ between 1 and $n-1$, let $D_{\frac {i}{n}}=D^i$, and we are done. 

We can apply a symmetric argument to produce a sequence of disks, $\{E_{\frac {i}{m}}\}_{i=1,...,m-1}$, such that $E_{\frac {i}{m}} \cap E_{\frac {i+1}{m}} =\emptyset$, for all $i$ between 0 and $m-1$.

Our goal now is to use all of these disks to define $\Phi$. The definition we give may seem overly technical, but there are certain features that are worth pointing out before-hand. If $\Phi$ is any map from $S \times D^2$ into $M$, then let $\Phi _x (s)=\Phi (s,x)$. If $x$ is the center of $D^2$, we would like $\Phi _x$ to be the identity on $S$. That is, at the center of $D^2$, we do nothing to $S$. For $x$ near the boundary of $D^2$, we would like $\Phi _x(S)$ to be disjoint from at least one of the disks, $\{D_{\frac {i}{n}}\}$, or $\{E_{\frac {i}{m}}\}$. In particular, there should be two points, $\theta _0$ and $\theta _1$ on $\partial D^2$ such that for all $x \in D$ near $\theta _i$, $\Phi _x (S)$ is disjoint from both $D_i$ and $E_i$, for $i=0,1$. In fact, any map $\Phi$ which fits this description will work for the remainder of our proof. We give an explicit construction of one here for concreteness. 

Let $\{U_i\}_{i=0,...,n}$ and $\{V_i\}_{i=0,...,m}$ denote neighborhoods of the disks, $\{D_{\frac {i}{n}}\}$ and $\{E_{\frac {i}{m}}\}$, such that
\begin{enumerate}
    \item $U_i \cap U_{i+1}=\emptyset$, for $0 \le i <n$
    \item $V_i \cap V_{i+1}=\emptyset$, for $0 \le i <m$
    \item $U_0 \cap V_0=\emptyset$
    \item $U_n \cap V_m=\emptyset$
\end{enumerate}

For each $i$ between 0 and $n$, let $\gamma ^i :M \times I \rightarrow M$ be an isotopy, such that $\gamma ^i _0 (x)=x$ for all $x \in M$, $\gamma ^i _t (x)=x$ for all $t$, and all $x$ outside of $U_i$, and $\gamma ^i_1(S) \cap D_{\frac {i}{n}} =\emptyset$ (where $\gamma ^i _t(x)$  is short-hand for $\gamma ^i (x,t)$). In other words, $\gamma ^i$ is an isotopy which pushes $S$ off of $D_{\frac {i}{n}}$, inside $U_i$. Similarly, for each $i$ between 0 and $m$, let $\delta ^i$ be an isotopy which pushes $S$ off of $E_{\frac {i}{m}}$, inside $V_i$.

Now choose $n+m+2$ points on $S^1$, which are labelled and cyclically ordered as follows: $x_0,...,x_n,y_m,...,y_0$. Let $\{f_i\}_{i=0,...,n}$ and $\{g_i\}_{i=0,...,m}$ be sets of continuous functions from $S^1$ to $I$, such that
\begin{enumerate}
    \item $supp(f_i)=(x_{i-1},x_{i+1})$, for $0<i<n$,
    \item $supp(g_i)=(y_{i+1},y_{i-1})$, for $0<i<m$,
    \item $supp(f_0)=(y_0,x_1)$, $supp(f_n)=(x_{n-1},y_m)$,
    \item $supp(g_0)=(y_1,x_0)$, $supp(g_m)=(x_n,y_{m-1})$,
    \item for each $\theta \in S^1$, there is an $i$ such that $f_i(\theta)=1$ or $g_i(\theta)=1$,
    \item there are points, $\theta _0 \in (y_0,x_0)$ and $\theta _1 \in (x_n,y_m)$, such that $f_0(\theta _0)=g_0(\theta _0)=f_n(\theta _1)=g_m(\theta _1)=1$.
\end{enumerate}
The existence of such functions is easily verified by the reader. 

Finally, we define $\Phi :S \times D^2 \rightarrow M$. For each $x \in S$ and $(r,\theta)\in D^2$, let 
\[\Phi _{(r,\theta)}(x)=\prod \limits _{i=0} ^{n} \gamma ^i _{rf_i(\theta)}(x) \prod \limits _{i=0} ^{m} \delta ^i _{rg_i(\theta)}(x)\]

Note that it is a consequence of how we defined the functions $\{f_i\}$ and $\{g_i\}$ that for each value of $\theta$ there are at most two functions that are not the identity. Since these two functions will have support on disjoint subsets of $M$, the order of the above product does not matter. For example, suppose $n=7$. Let $\theta '$ be some value of $S^1$ between $x_3$ and $x_4$. Then $f_3$ and $f_4$ are the only functions that can be non-zero at $\theta '$. Hence, the above product simplifies to:
\[\Phi _{(r,\theta')}(x)=\gamma ^3 _{rf_3(\theta ')}(x) \cdot \gamma ^4 _{rf_4(\theta ')}(x)\]
But $\gamma ^3$ is the identity outside $U_3$, $\gamma ^4$ is the identity outside $U_4$, and $U_3 \cap U_4 =\emptyset$. Hence, $\gamma ^3$ and $\gamma ^4$ commute. 

We now perturb $\Phi$ slightly so that it is in general position with respect to $S$, and again denote the new function as $\Phi$. Consider the set $\Sigma=\{x \in D^2 | \Phi _x (S)$ is not transverse to $F\}$. If $\Phi$ is in general position with respect to $S$, then Cerf theory (see \cite{cerf:68}) tells us that $\Sigma$ is homeomorphic to a graph, and the maximum valence of each vertex of this graph is 4. We will use these facts later.

\medskip
\noindent \underline{Stage 2:} {\it Labelling $D^2$}.
\medskip

Let $C_0$ and $C_1$ denote sets of components of $\Gamma (F)$ such that 
\begin{enumerate}
    \item any isolated vertex is in both sets,
    \item a component which is not an isolated vertex (i.e. it contains at least one edge) is in exactly one of the sets, and
    \item the component containing the edge $D_i-E_i$ is an element of $C_i$, for $i=0,1$.
\end{enumerate}

\medskip

\noindent {\bf Note:} The salient feature of these sets is that if $v$ is a vertex of $\Gamma (F)$ in a component, $A \in C_0$, and $v'$ is some other vertex, in a component, $B \in C_1$, then $A \ne B$. Hence, $v-v'$ cannot be an edge of $\Gamma (F)$. We may conclude then that if $D$ is a red disk representing $v$, and $E$ is a blue disk representing $v'$, then $D \cap E \ne \emptyset$. We will make extensive use of these facts. 
\medskip

As at the end of Stage 1, let $\Sigma=\{x \in D^2 | \Phi _x (S)$ is not transverse to $F\}$. A {\it region} of $D^2$ is a component of $D^2 - \Sigma$. Let $x$ be any point in the interior of some region. Label this region with a ``$D_i$" (``$E_i$") if there is a disk, $\Delta$, in $\Phi _x (S)$ such that $\partial \Delta$ is an essential loop on $F$, $int (\Delta) \cap F$ is a (possibly empty) union of inessential loops on $F$, and $\Delta$ is isotopic rel $\partial$ to a red (blue) compressing disk for $F$ which corresponds to a vertex of $\Gamma (F)$ in $C_i$. It is possible that some regions will have more than one label. For example, if there is a red disk in $\Phi _x (S)$ which corresponds to an isolated vertex of $\Gamma (F)$, then it will have both of the labels ``$D_0$" and ``$D_1$".

\medskip
\noindent {\bf Note:} We will be consistent about using quotation marks to denote labels. Hence, the label ``$D_0$", for example, should not be confused with the disk, $D_0$.
\medskip

The point of our labelling is that if there is any compressing disk for $F$ contained in $\Phi _x (S)$, then the region containing $x$ will be assigned a label. Hence, the existence of an unlabelled region would imply that all components of intersection of $\Phi _x (S)$ and $F$ are either essential on both surfaces, or inessential on both. From there it is easy to prove the Theorem. We now proceed under the assumption that there are no unlabelled regions.

\begin{clm}
\label{c:same}
No region can have both of the labels ``$D_i$" and ``$E_{1-i}$".
\end{clm}

\begin{proof}
Let $x$ be a point in the interior of a region with the labels ``$D_0$" and ``$E_1$". Then there is a red disk, $D$, and a blue disk, $E$, in $\Phi _x (S)$. Furthermore, according to our labelling, the equivalence class of $D$ represents a vertex of $\Gamma (F)$ in a component $A \in C_0$. Similarly, $[E] \in B$, for some $B \in C_1$. But since $D$ and $E$ are both subsets of $\Phi _x (S)$ they must be disjoint, a contradiction (see the above note after the definitions of $C_0$ and $C_1$). 
\end{proof}

\begin{clm}
\label{c:adjacent}
If a region has the label ``$D_i$", then no adjacent region can have the label ``$E_{1-i}$".
\end{clm}

{\bf Note.} The proof of this claim is essentially the same as that of Lemma 4.4 of \cite{gabai:87}. 

\begin{proof}
Suppose the region, $\mathcal R_0$, has the label ``$D_0$", and an adjacent region, $\mathcal R_1$, has the label ``$E_1$". For each $x \in \mathcal R_0$, there is a red disk, $D$, contained in $\Phi _x (S)$. Similarly, for each $x \in \mathcal R_1$, there is a blue disk, $E$, contained in $\Phi _x (S)$.

Let $x_i$ be some point in the interior of $\mathcal R_i$. Let $p:I \rightarrow D^2$ be an embedded path connecting $x_0$ to $x_1$, which does not wander into any region other than $\mathcal R_0$ or $\mathcal R_1$. As $t$ increases from $0$ to $1$, we see a moment, $t_*$, when $\Phi _{p(t_*)}(S)$ does not meet $F$ transversely (i.e. $t_*=p^{-1}(\partial \mathcal R_1)$). At $t_*$ we simultaneously see the disappearance of $D$, and the appearance of $E$. (Otherwise, $\mathcal R_1$ would have both the labels ``$D_0$" and ``$E_1$", which we ruled out in Claim \ref{c:same}).

Let $\gamma _t$ be the components of $\Phi _{p(t)}(S) \cap F$, whose intersection with $D$ is nontrivial if $t \le t_*$, and whose intersection with $E$ is nontrivial for $t>t_*$. As $t$ approaches $t_*$ from either side, we see a tangency occur for $\gamma _t$ (either with itself, or with some other component of $\Phi _{p(t)}(S) \cap F$). Since generically only one such tangency occurs at $t_*$, we see that $\lim \limits _{t \rightarrow t_*^-} \gamma_t \cap \lim \limits _{t \rightarrow t_*^+} \gamma_t \ne \emptyset$, as in Figure \ref{f:gamma}. 

        \begin{figure}[htbp]
        \psfrag{g}{$\gamma _t$}
        \psfrag{-}{$t<t_*$}
        \psfrag{+}{$t>t_*$}
        \psfrag{0}{$t=t_*$}
        \vspace{0 in}
        \begin{center}
        \epsfxsize=3 in
        \epsfbox{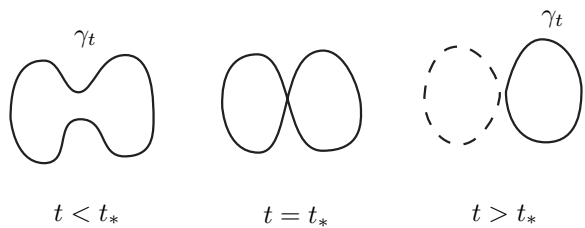}
        \caption{$\gamma _t$, before and after $t_*$.}
        \label{f:gamma}
        \end{center}
        \end{figure}

Since $D$ and $E$ are on opposite sides of $F$, we see from Figure \ref{f:gamma} that they can be made disjoint (since $F$ is orientable), and hence, $D - E$ is an edge of $\Gamma (F)$. As in the proof of Claim \ref{c:same}, this contradicts our labelling. 
\end{proof}

\medskip
\noindent \underline{Stage 3:} {\it The 2-complex, $\Pi$, and a map from $D^2$ to $\Pi$}.
\medskip
                                                                        
Let $\Pi$ be the labelled 2-complex depicted in Figure \ref{f:pi}. Let $\Sigma '$ be the dual graph of $\Sigma$. Map each vertex of $\Sigma '$ to the point of $\Pi$ with the same label(s) as the region of $D^2$ in which it sits. Claim \ref{c:same} assures that this map is well defined on the vertices of $\Sigma '$. 

        \begin{figure}[htbp]
        \psfrag{1}{``$D_0$"}
        \psfrag{2}{``$D_0$",``$D_1$"}
        \psfrag{3}{``$D_1$"}
        \psfrag{4}{``$E_0$"}
        \psfrag{5}{``$E_0$",``$E_1$"}
        \psfrag{6}{``$E_1$"}
        \psfrag{7}{``$D_0$",``$E_0$"}
        \psfrag{8}{``$D_1$",``$E_1$"}
        \vspace{0 in}
        \begin{center}
        \hspace{-.75in}
        \epsfxsize=4 in
        \epsfbox{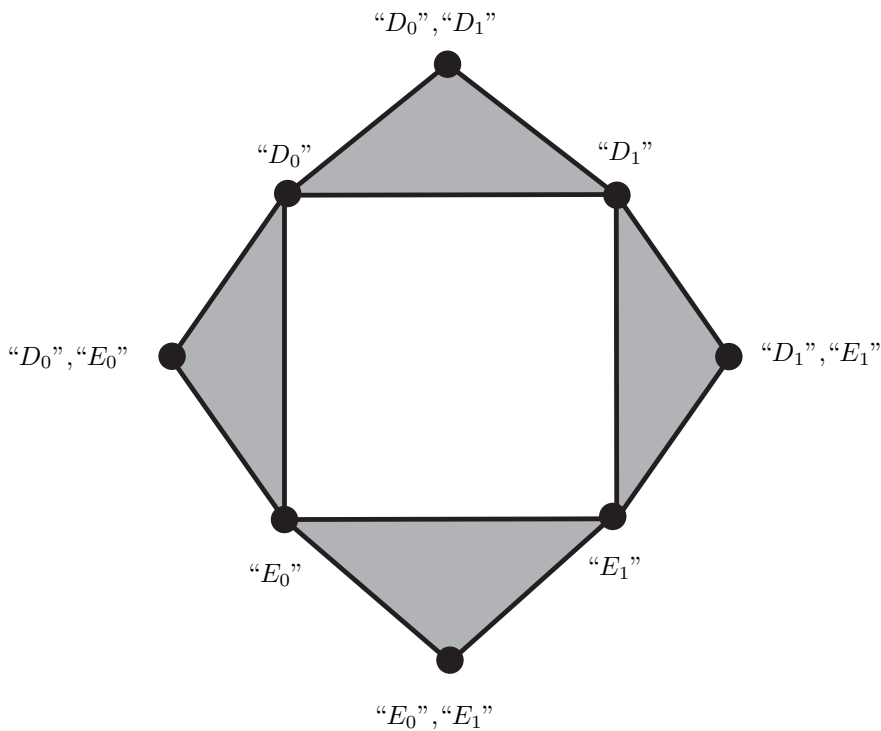}
        \caption{The 2-complex, $\Pi$.}
        \label{f:pi}
        \end{center}
        \end{figure}

Similarly, map each edge of $\Sigma '$ to the 1-simplex of $\Pi$ whose endpoints are labelled the same. Claim \ref{c:adjacent} guarantees that this, too, is well defined. 

We now claim that the map to $\Pi$ extends to all of $D^2$. Note that the maximum valence of a vertex of $\Sigma$ is four. Hence, the boundary of each region in the complement of $\Sigma '$ gets mapped to a 1-cycle with at most four vertices in $\Pi$. Inspection of Figure \ref{f:pi} shows that there is only one such cycle which is not null homologous. Hence, to show the map extends, it suffices to prove the following:

\begin{clm}
All four labels cannot occur around a vertex of $\Sigma$.
\end{clm}

\begin{proof}
The proof is somewhat similar to that of Claim \ref{c:adjacent}. Let $x_*$ denote a vertex of $\Sigma$, around which we see all four labels. Claims \ref{c:same} and \ref{c:adjacent} imply that each region around such a vertex must have a unique label. If $x$ is any point in a region labelled ``$D_i$" (for $i=0,1$), then there is a red disk, $D'_i$, contained in $\Phi _x (S)$. Similarly, if $x$ is in a region labelled ``$E_i$", then there is a blue disk, $E'_i$, contained in $\Phi _x (S)$.

Let $\gamma _x$ denote the components of $\Phi _x (S) \cap F$ whose intersection with $D'_i$ is nontrivial for $x$ in a region with the label ``$D_i$", and whose intersection with $E'_i$ is nontrivial for $x$ in a region with the label ``$E_i$". Let $\gamma _{x_*}$ denote the union of all possible limits of $\gamma _x$, as $x$ approaches $x_*$ from different directions. Then $\gamma _{x_*}$ is some graph with exactly 2 vertices, each of valence four. If $\mathcal R$ is a region which meets $x_*$, and $x \in \mathcal R$, then $\gamma _x$ is one component of the 1-manifold obtained from $\gamma _{x_*}$ by resolving both vertices in some way. (To resolve a vertex is to replace a part of $\gamma _{x_*}$ that looks like ``$\times$" by a part that looks like ``$)($".)

Since there are exactly four ways to resolve two vertices, and there are four regions with different labels around $x_*$, we see all four resolutions of $\gamma _{x_*}$ appear around $x_*$. However, the orientability of $F$ insures that {\it some} resolution of $\gamma _{x_*}$ always produces a 1-manifold, $\alpha$, that can be made disjoint from all components of all other resolutions. (In fact, one can show that there are always two resolutions of $\gamma _{x_*}$ which have this property, but this is irrelevant for us). Some component of $\alpha$ is the boundary of one of the disks, $D'_0$, $D'_1$, $E'_0$, or $E'_1$. Suppose $\partial D'_0 \subset \alpha$. Then $\partial D'_0$ can be made disjoint from $\partial E'_1$, since $\partial E'_1$ is a component of some other resolution of $\gamma _{x_*}$. This now contradicts our labelling. 
\end{proof}

\medskip
\noindent \underline{Stage 4:} {\it Finding an unlabelled subregion}.
\medskip

To obtain a contradiction, it suffices to prove that the map from $D^2$ to $\Pi$, when restricted to $\partial D^2$, induces a non-trivial map on homology. To this end, we must examine the possibilities for the labels of the regions adjacent to $\partial D^2$. 

Recall from Stage 1 the points, $\{x_i\}$ and $\{y_i\}$ on $S^1=\partial D^2$, and the disks, $\{D _{\frac{i}{n}}\}$ and $\{E _{\frac{i}{m}}\}$. 

\begin{clm}
\label{c:boundary1}
If $\mathcal R$ is a region such that some point, $p$, of $\partial \mathcal R$ lies between $x_i$ and $x_{i+1}$ on $\partial D^2$ (for some $i$), then $\mathcal R$ cannot have both of the labels ``$E_0$" and ``$E_1$".
\end{clm}

\begin{proof}
As discussed in Stage 1, at $p$ the definition of $\Phi$ simplifies to 
\[\Phi _p(x)=\gamma ^i _{f_i(p)}(x) \cdot \gamma ^{i+1} _{f_{i+1}(p)}(x)\]

By our construction of the functions, $\{f_i\}$, we know that either $f_i (p)=1$ or $f_{i+1}(p)=1$. Suppose the former is true. Then $\Phi$ further simplifies to
\[\Phi _p(x)=\gamma ^i _1(x) \cdot \gamma ^{i+1} _{f_{i+1}(p)}(x)\]

By construction, $\gamma ^i _1(S) \cap D_{\frac {i}{n}}=\emptyset$, and hence $\Phi _p (S) \cap D_{\frac {i}{n}}=\emptyset$. By our definition of a {\it region}, we conclude that for every point, $p \in \mathcal R$, $\Phi _p (S) \cap D_{\frac {i}{n}}=\emptyset$.

Now, if $\mathcal R$ has both of the labels ``$E_0$" and ``$E_1$", then either
\begin{enumerate}
    \item there are blue disks, $E$ and $E'$, in $\Phi _p (S)$, which correspond to vertices in different components of $\Gamma (F)$, or
    \item there is a blue disk, $E \subset \Phi _p (S)$, which corresponds to an isolated vertex of $\Gamma (F)$.
\end{enumerate}
In either case, there is a blue disk (say $E$) which corresponds a vertex of $\Gamma (F)$ which is in a component other than the one that contains $D_{\frac {i}{n}}$. However, $\Phi _p (S) \cap D_{\frac {i}{n}}=\emptyset$ and $E \subset \Phi _p (S)$ implies that $D_{\frac {i}{n}}-E$ is an edge of $\Gamma (F)$, a contradiction. 
\end{proof}

By a symmetric argument, we can show that a region adjacent to a point of $\partial D^2$ between $y_0$ and $y_m$ cannot have both of the labels ``$D_0$" and ``$D_1$".

\begin{clm}
\label{c:boundary2}
If $\mathcal R_0$ and $\mathcal R_1$ are regions adjacent to $\partial D^2$, with labels ``$E_0$" and ``$E_1$", respectively, then no point of $\mathcal R_0 \cap \mathcal R_1$ can lie between $x_i$ and $x_{i+1}$ on $\partial D^2$ (for any $i$).
\end{clm}

\begin{proof}
If the Claim is not true, then some point, $p$, of $\mathcal R_0 \cap \mathcal R_1$ lies between $x_i$ and $x_{i+1}$ on $\partial D^2$, for some $i$. As in the proof of Claim \ref{c:boundary1}, we may conclude that $\Phi _p (S)$ is disjoint from the red compressing disk, $D_{\frac {i}{n}}$. However, this implies that for all $x \in D^2$ near $p$, $\Phi _x (S) \cap D_{\frac {i}{n}}=\emptyset$. In particular, there are points, $p_j \in \mathcal R_j$ (for $j=0,1$), such that $\Phi _{p_j} (S) \cap D_{\frac {i}{n}}=\emptyset$.

Since $\mathcal R_j$ has only the label ``$E_j"$ (for $j=0,1$), there are blue disks, $E'_j \subset \Phi _{p_j}(S)$, which are in different components of $\Gamma (F)$. But $E'_0 \cap D_{\frac {i}{n}}=E'_1 \cap D_{\frac {i}{n}}=\emptyset$ implies that $E'_0 - D_{\frac {i}{n}} - E'_1$ is a chain of edges in $\Gamma (F)$, a contradiction. 
\end{proof}

Again by a symmetric argument, we can show that there cannot be a pair of regions with the labels ``$D_0$" and ``$D_1$", whose boundaries both contain a point of $\partial D^2$ between $y_0$ and $y_m$.

Now, recall from our definition of the functions $\{f_i\}$ and $\{g_i\}$ that there is a  point, $\theta _0 \in (y_0,x_0)$, such that $f_0(\theta _0)=g_0(\theta _0)=1$. 

\begin{clm}
\label{c:boundary3}
If $\mathcal R_0$ is the region whose boundary contains $\theta _0$, then $\mathcal R_0$ cannot have either of the labels ``$D_1$" or ``$E_1$". 
\end{clm}

\begin{proof}
At $\theta _0$, the definition of $\Phi$ simplifies to
\[\Phi _{\theta _0}(x)=\gamma ^0 _1(x) \cdot \delta ^0 _1(x)\]

Hence, $\Phi _{\theta _0}(S)$ is disjoint from both $D_0$ and $E_0$. If $\mathcal R_0$ had the label ``$D_1$", then there would be some red disk, $D \subset \Phi _{\theta _0}(S)$, which corresponds to a vertex of $\Gamma (F)$ in some component, $B \in C_1$. Any such disk must intersect $E_0$, since $E_0$ corresponds to a vertex of $\Gamma (F)$ in some component, $A \in C_0$. However, $D \cap E_0 =\emptyset$, a contradiction. A symmetric argument rules out the possibility of the label ``$E_1$" for $\mathcal R_0$. 
\end{proof}

Similarly, we can show that if $\mathcal R_1$ is the region whose boundary contains $\theta _1$, then $\mathcal R_1$ cannot have either of the labels ``$D_0$" or ``$E_0$".

\begin{clm}
The map from $D^2$ to $\Pi$, when restricted to $\partial D^2$, induces a non-trivial map on homology.
\end{clm}

\begin{proof}
Claims \ref{c:boundary1} and \ref{c:boundary2} imply that no point of the arc of $\partial D^2$ between $x_0$ and $x_n$ gets mapped to the lower triangle of $\Pi$. Similarly, we can show that no point of the arc of $\partial D^2$ between $y_0$ and $y_m$ gets mapped to the upper triangle. Finally, Claim \ref{c:boundary3} implies that the region, $\mathcal R_0$, whose boundary contains the point, $\theta _0$, gets mapped to the left triangle of $\Pi$, while the region, $\mathcal R_1$, whose boundary contains the point, $\theta _1$, gets mapped to the right triangle. 

All of this directly implies the Claim, as illustrated in Figure \ref{f:rpi}.

        \begin{figure}[htbp]
        \psfrag{1}{$y_0$}
        \psfrag{2}{$x_0$}
        \psfrag{3}{$x_n$}
        \psfrag{4}{$y_m$}
        \vspace{0 in}
        \begin{center}
        \epsfxsize=3 in
        \epsfbox{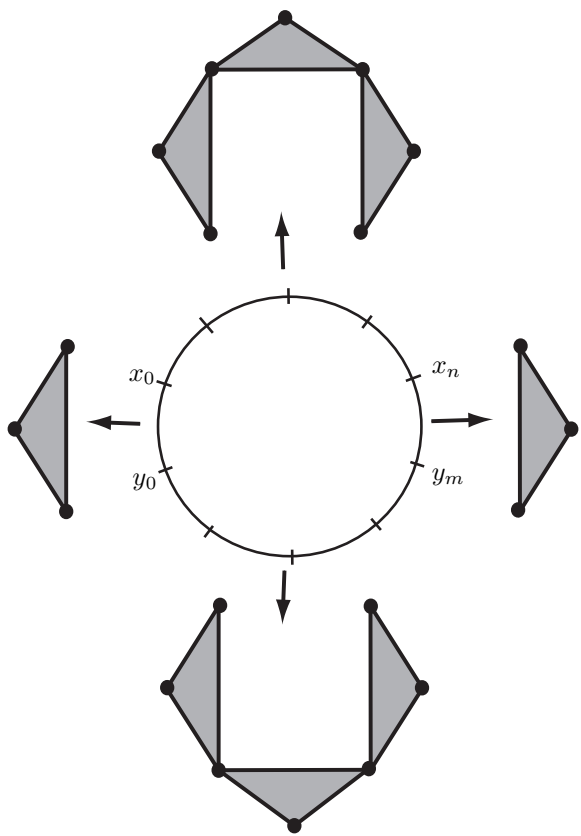}
        \caption{The map from $\partial D^2$ to $\Pi$ is non-trivial on homology.}
        \label{f:rpi}
        \end{center}
        \end{figure}

\end{proof}

We now complete the proof of Theorem \ref{t:essential_intersection}. Since $H_1(\Pi) \ne 0$, the map to $\Pi$ cannot extend, so there must be a region with no labels. Now, let $x$ be a point in the interior of an unlabelled region. Then there are no compressions in $\Phi _x(S)$ for $F$. Hence, any loop of $F \cap \Phi _x(S)$ is either essential on both surfaces, or inessential on both. It is now a routine matter to use an innermost disk argument to show that the inessential loops can be removed by a sequence of disk swaps, or by a further isotopy of $S$, if $M$ is irreducible. 
\end{proof}

We now present some immediate corollaries to Theorem \ref{t:essential_intersection}. 

\begin{cor}
If $M_1 \# M_2$ contains a critical surface, then either $M_1$ or $M_2$ contains a critical surface. 
\end{cor}

\begin{cor}
A reducible 3-manifold does not admit a critical Heegaard splitting.
\end{cor}

\begin{cor}
\label{c:essential_boundary}
Suppose $M$ is a 3-manifold which admits a critical Heegaard splitting, such that $\partial M \ne \emptyset$. Then $\partial M$ is essential in $M$. 
\end{cor}

\begin{proof}
If $\partial M \cong S^2$, then either $\partial M$ is essential, or $M \cong B^3$. However,  Theorem \ref{t:nocrit} implies that $B^3$ does not contain a critical Heegaard surface. 

If $\partial M$ contains a component of non-zero genus, then we claim it is incompressible. If not, then Theorem \ref{t:essential_intersection} implies that there is a compressing disk for $\partial M$ which misses the critical Heegaard splitting. However, this would imply that there is a compression body, $W$, with $\partial _- W$ compressible in $W$, a contradiction. 
\end{proof}

\section{Generalized Heegaard Splittings}

In this section we cover the background material that we will need for the proof of our main result, Theorem \ref{t:common_stab}.

\begin{dfn} 
(Scharlemann- Thompson \cite{st:94}) 
A {\it Generalized Heegaard Splitting} (GHS) of a 3-manifold, $M$, is a sequence of closed, embedded, pairwise disjoint surfaces, $\{F_i\}_{i=0}^{2n}$, such that for each $i$ between 1 and $n$, $F_{2i-1}$ is a non-trivial Heegaard splitting, or a union of Heegaard splittings (at least one of which is non-trivial), of the submanifold of $M$ co-bounded by $F_{2i-2}$ and $F_{2i}$, and such that $\partial M=F_0 \amalg F_{2n}$. 
\end{dfn}

{\it Notes:} (1) If $M$ is closed, then $F_0=F_{2n}=\emptyset$. (2) We allow $F_{2i-1}$ to be a union of Heegaard splittings only when the submanifold of $M$ co-bounded by $F_{2i-2}$ and $F_{2i}$ is disconnected. 

\begin{figure}[htbp]
        \psfrag{1}{$F_0$}
        \psfrag{2}{$F_1$}
        \psfrag{3}{$F_2$}
        \psfrag{4}{$F_{2n-1}$}
        \psfrag{5}{$F_{2n}$}
        \vspace{0 in}
        \begin{center}
        \epsfxsize=1 in
        \epsfbox{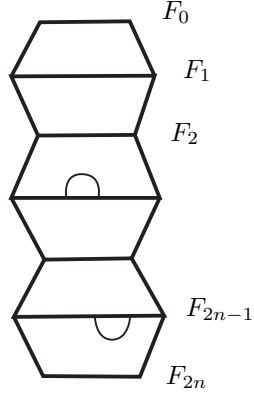}
        \caption{Schematic depicting a Generalized Heegaard Splitting.}
        \label{f:ghs}
        \end{center}
        \end{figure}

We will sometimes depict a GHS schematically as in Figure \ref{f:ghs}. Often when we do this we will also need to represent compressing disks for $F_{2i-2}$, for various values of $i$. Examples of this are the curved arcs depicted in the figure.

Recall the ordering of embedded surfaces given in Definition \ref{d:order}. We now use this to define a partial ordering on GHSs. 

\begin{dfn}
Let $F^1=\{F^1_i\}$ and $F^2=\{F^2_j\}$ be two GHSs of a 3-manifold, $M$. We say $F^1<F^2$ if $\{F^1_i | i \ {\rm odd}\} <\{F^2_j | j \ \rm odd\}$, where each set is put in non-increasing order, and then the comparison is made lexicographically. 
\end{dfn}

We now define two ways to get from one GHS to a smaller one. Suppose $\{F_i\}$ is a GHS. Suppose further that $D-E$ is an edge in $\Gamma (F_{2i-1})$, where $D$ and $E$ are disks in the submanifold of $M$ co-bounded by $F_{2i-2}$ and $F_{2i}$, for some $i$. Let $F_D$ denote the surface obtained from $F_{2i-1}$ by compression along $D$, and $F_E$ denote the surface obtained from $F_{2i-1}$ by compression along $E$. If $D \cap E =\emptyset$, then let $F_{DE}$ denote the surface obtained from $F_{2i-1}$ by compression along both $D$ and $E$. There are now two cases, with several subcases:

\begin{enumerate}
    \item $D \cap E =\emptyset$
        \begin{enumerate}
            \item $F_D \ne F_{2i-2}$, $F_E \ne F_{2i}$. 
            \\ Remove $F_{2i-1}$ from $\{F_i\}$. In it's place, insert $\{F_D, F_{DE}, F_E\}$ and reindex.  
            \item $F_D = F_{2i-2}$, $F_E \ne F_{2i}$. 
            \\ Replace $\{F_{2i-2}, F_{2i-1}\}$ with $\{F_{DE},F_E\}$.
            \item $F_D \ne F_{2i-2}$, $F_E = F_{2i}$. 
            \\ Replace $\{F_{2i-1}, F_{2i}\}$ with $\{F_D,F_{DE}\}$.
            \item $F_D = F_{2i-2}$, $F_E = F_{2i}$. 
            \\ Replace $\{F_{2i-2}, F_{2i-1}, F_{2i}\}$ with $F_{DE}$ and reindex. 
        \end{enumerate}
    \item $|D \cap E|=1$ (In this case $F_D$ and $F_E$ co-bound a product region of $M$)
        \begin{enumerate}
            \item $F_D \ne F_{2i-2}$, $F_E \ne F_{2i}$.
            \\ Replace $F_{2i-1}$ in $\{F_i\}$ with $F_D$ or $F_E$.
            \item $F_D = F_{2i-2}$, $F_E \ne F_{2i}$. 
            \\ Remove $\{F_{2i-2}, F_{2i-1}\}$ and reindex. 
            \item $F_D \ne F_{2i-2}$, $F_E = F_{2i}$.
            \\ Remove $\{F_{2i-1}, F_{2i}\}$ and reindex. 
            \item $F_D = F_{2i-2}$, $F_E = F_{2i}$.
            \\ Remove $\{F_{2i-1}, F_{2i}\}$ {\it or} $\{F_{2i-2}, F_{2i-1}\}$ and reindex. 
        \end{enumerate}
\end{enumerate}

If $\partial D \cup \partial E$ bounds an annulus on $F_{2i-1}$, then some element of the sequence of surfaces that we get after one of the above operations contains a 2-sphere component. In this case, we remove the 2-sphere component. 

\begin{figure}[htbp]
        \psfrag{a}{(a)}
        \psfrag{b}{(b)}
        \vspace{0 in}
        \begin{center}
        \epsfxsize=4 in
        \epsfbox{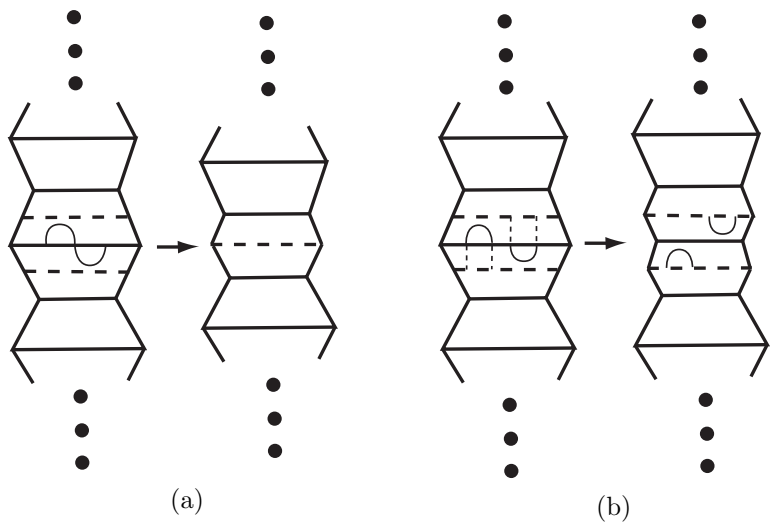}
        \caption{(a) A destabilization. (b) A weak reduction.}
        \label{f:reddfn}
        \end{center}
        \end{figure}

In Case 1 above ($D \cap E =\emptyset$), we say the new sequence was obtained from the old one by the {\it weak reduction}, $D-E$. In Case 2 ($|D \cap E|=1$), we say the new sequence was obtained by a {\it destabilization}. Each of these operations is represented schematically in Figure \ref{f:reddfn}. In either case we leave it as an exercise to show that the new sequence is a GHS, provided $M$ is irreducible. Note that if the GHS, $F^1$, is obtained from the GHS, $F^2$, by weak reduction or destabilization, then $F^1<F^2$. 

For readers unfamiliar with Generalized Heegaard Splittings, we pause here for a moment to tie these concepts to more familiar ones. Suppose $M$ is a closed 3-manifold, and $\{F_i\}_{i=0}^{2n}$ is a GHS of $M$. Since $M$ is closed, we must have $F_0=F_{2n}=\emptyset$. If, in addition, $n=1$, then our GHS looks like $\{\emptyset, F_1, \emptyset\}$. By definition, $F_1$ is a Heegaard splitting of $M$. If $\{\emptyset, F'_1, \emptyset\}$ was obtained from $\{\emptyset, F_1, \emptyset\}$ by a destabilization, then the Heegaard splitting, $F_1$, is a stabilization of the Heegaard splitting, $F'_1$.

\section{Minimal Genus Common Stabilizations}

In this section we prove the main result of this paper:

\begin{thm}
\label{t:common_stab}
Suppose $M$ is a 3-manifold whose boundary, if non-empty, is incompressible, and $F$ and $F'$ are distinct strongly irreducible Heegaard splittings of $M$ which induce the same partition of $\partial M$. If the minimal genus common stabilization of $F$ and $F'$ is not critical, then $M$ contains a non-boundary parallel incompressible surface.
\end{thm}

We will actually prove the contrapositive of this Theorem. That is, we show that if $F$ and $F'$ are strongly irreducible Heegaard splittings of a small 3-manifold, $M$, then their minimal genus common stabilization is critical. A {\it small} 3-manifold is one which is either closed and non-Haken, or is one with incompressible boundary, in which every incompressible surface is boundary parallel. 

They key technique that will be used in the proof is a careful analysis of the ``best" way one can transform $F$ into $F'$ by a sequence of intermediate GHSs. 

\begin{dfn}
A {\it Sequence Of GHSs} (SOG) of a 3-manifold, $M$, is a sequence, $\{F^j\}_{j=1}^n$, such that for each $k$ between 1 and $n-1$, one of the GHSs, $F^k$ or $F^{k+1}$, is obtained from the other by a weak reduction or destabilization.
\end{dfn}

{\it Notation:} We will always use subscripts to denote surfaces, superscripts to denote GHSs, and a boldface font to denote an entire SOG. Hence, $F_i^j$ is the $i$th surface of the $j$th GHS of the SOG, $\bf F$. If $\bf F$ is a SOG of a 3-manifold, $M$, then $M^k_i$ will always denote the submanifold of $M$ cobounded by $F^k_{i-1}$ and $F^k_{i+1}$.

\begin{dfn}
If $\bf F$ is a SOG of $M$, and $k$ is such that $F^{k-1}$ and $F^{k+1}$ are both obtained from $F^k$ by weak reduction or destabilization, then we say the GHS, $F^k$, is {\it maximal} in $\bf F$. Similarly, if $k$ is such that $F^k$ is obtained from both $F^{k-1}$ and $F^{k+1}$ by weak reduction or destabilization, then we say it is {\it minimal}.
\end{dfn}

We now define a partial ordering on SOGs. 

\begin{dfn}
Let $\bf F$ and $\bf G$ be two SOGs of some 3-manifold, $M$. Let $MAX(\bf F)$ and $MAX(\bf G)$ denote the sets of maximal GHSs of $\bf F$ and $\bf G$. If we reorder each of these sets in non-increasing order, then we say $\bf F<\bf G$ if $MAX({\bf F})<MAX({\bf G})$, where all comparisons are made lexicographically. If $\bf F$ is minimal among all SOGs of some collection, then we say it is {\it flattened}. 
\end{dfn}

This definition is a bit hard to slog through. After all, $MAX(\bf F)$ is a set of sets of surfaces. At every level we make a lexicographic comparison. Unfortunately, there does not seem to be much that can be done to make this any simpler. Perhaps the only thing that will convince the reader that this is the right definition is the proofs of a few Lemmas.

Let us not forget that our eventual goal is to prove Theorem \ref{t:common_stab}. To this end, we suppose for the remainder of the paper that $M$ is a small 3-manifold that contains non-isotopic strongly irreducible Heegaard splittings, $F$ and $F'$. One subtle technical point here is that if $M$ is not closed, then we will assume further that $F$ and $F'$ induce the same partition of $\partial M$. That is, we assume that $\partial M=\partial _1 M \cup \partial _2 M$ (either of which may be empty), and that all components of $\partial _1 M$ are on one side of both $F$ and $F'$, and all components of $\partial _2 M$ are on the other side of $F$ and $F'$. 

Among all SOGs of $M$ whose first GHS is $\{\partial _1 M, F, \partial _2 M \}$, and last GHSs is $\{\partial _1 M, F', \partial _2 M \}$, let $\bf F$ denote one which is flattened. We now prove several Lemmas that tell us precisely what $\bf F$ looks like.

\begin{lem}
\label{l:oddstrirr}
If $F^k$ is maximal in $\bf F$, then for every odd value of $i$ except for exactly one, $F^k_i$ is strongly irreducible in $M^k_i$.
\end{lem}

\begin{proof}
Since $F^k$ is maximal in $\bf F$, there is some odd $p$ such that $F^{k-1}$ is obtained from $F^k$ by a weak reduction or destabilization corresponding to an edge, $D-E$, in $\Gamma (F^k_p)$. Similarly, there is an odd $q$ such that $F^{k+1}$ is obtained from $F^k$ by a weak reduction or destabilization corresponding to an edge, $D'-E'$, in $\Gamma (F^k_q)$. We first claim that $p=q$. If not, then we can replace $F^k$ in ${\bf F}$ with the GHS obtained from $F^{k-1}$ by the weak reduction or destabilization, $D'-E'$. Since the new GHS can also be obtained from $F^{k+1}$ by the weak reduction or destabilization, $D-E$, our substitution has defined a new SOG of $M$. In general, this replaces one maximal GHS of ${\bf F}$ with two smaller maximal GHSs, reducing $MAX({\bf F})$, and therefore contradicting our assumption that ${\bf F}$ was flattened. 

Now, suppose $F^k_r$ is not strongly irreducible, for some odd $r \ne p$, and let $D^*-E^*$ be an edge in $\Gamma (F^k_r)$. Let $G^-$, $G^0$, and $G^+$ denote the GHSs obtained from $F^{k-1}$, $F^k$, and $F^{k+1}$ by the weak reduction or destabilization, $D^*-E^*$. Now replace $F^k$ in ${\bf F}$ with $\{G^-,G^0,G^+\}$ and reindex. In general, this also replaces one maximal GHS of ${\bf F}$ with two smaller maximal GHSs, reducing $MAX({\bf F})$. 
\end{proof}

\begin{lem}
\label{l:critmax}
If $F^k$ is maximal in ${\bf F}$, then there is exactly one odd number, $p$, such that $F^k_p$ is critical in $M^k_p$.
\end{lem}

\begin{proof}
Let $p$, $D-E$, and $D'-E'$ be as defined in the proof of Lemma \ref{l:oddstrirr}. To establish the Lemma we show that $D-E$ and $D'-E'$ are in different components of $\Gamma (F^k_p)$. Suppose this is not true. Then there is a chain of edges in $\Gamma (F^k_p)$ connecting $D-E$ to $D'-E'$. That is, there is a sequence of disks, $\{C_l\}_{l=1}^w$, such that $C_1=D$, $C_2=E$, $C_{w-1}=D'$, and $C_w=E'$, and for each $l$, $[C_{l-1}]-[C_l]$ is an edge of $\Gamma (F^k_p)$ ($[C_l]$ denotes the equivalence class of $C_l$). In other words, the following is a chain connecting $D-E$ to $D'-E'$:
\[ [D]-[E]-[C_3]-[C_4]-...-[C_{w-2}]-[D']-[E']\]

\begin{clm}
\label{c:chain}
There are disks, $\{D_{2l+1}\}$ and $\{E_{2l}\}$, such that the following is a chain in $\Gamma (F^k_p)$: 
\[D_1-E_2-D_3-E_4-...-D_{w-1}-E_w\]
where $D_1=D$, $E_2=E$, and the pair $(D_{w-1},E_w)$ is isotopic to the pair $(D',E')$. 
\end{clm}

\begin{proof}
Let $D_1=D$, and $E_2=E$. By definition there exist representatives, $U \in [E_2]$, and $V \in [C_3]$ such that $|U \cap V|\le 1$. Since $E_2$ and $U$ are in $[E_2]$, there is an isotopy, $\Phi$, of $M$ taking $F^k_p$ to $F^k_p$, and $E_2$ to $U$. If we apply $\Phi ^{-1}$ to $V$, we obtain a disk, on the opposite side of $F^k_p$ as $E_2$, which intersects $E_2$ at most once. Let this disk be $D_3$. Then $E_2-D_3$ is an edge of $\Gamma (F^k_p)$, and $D_3 \in [C_3]$. Continue in this way, using $D_3$ now to construct $E_4$, etc. 
\end{proof}

\begin{clm}
\label{c:isotopy}
The GHS, $F'$, obtained from $F^k$ by the weak reduction or destabilization, $D_{w-1}-E_w$, is the same as $F^{k+1}$ (up to isotopy). 
\end{clm}

\begin{proof}
Let $\Phi$ denote the isotopy of $M$ which takes the pair, $(D_{w-1},E_w)$, to the pair, $(D',E')$. Then $\Phi$ also takes $F'$ to $F^{k+1}$.
\end{proof}

\begin{clm}
\label{c:disjoint}
For each even $l$, we may assume that $D_{l-1} \cap D_{l+1}=\emptyset$. Similarly, for each odd $l$, we may assume that $E_{l-1} \cap E_{l+1}=\emptyset$.
\end{clm}

\begin{proof}
If the statement of the claim is not true, then we show that we can replace the chain with a chain where it is true. Let us assume then that for some even $l$, $D_{l-1} \cap D_{l+1} \ne \emptyset$. 

Begin by isotoping $D_{l-1}$ and $D_{l+1}$ so that $|D_{l-1} \cap D_{l+1}|$ is minimal. Let $\alpha$ be an outermost arc of $D_{l-1} \cap D_{l+1}$ on $D_{l+1}$. Then $\alpha$ cuts off a subdisk, $V$, of $D_{l+1}$, such that $V \cap D_{l-1} = \alpha$. Furthermore, $\alpha$ cuts $D_{l-1}$ into two disks, $V'$ and $V''$. Since $D_{l-1}$ and $D_{l+1}$ each meet $E_l$ in at most a point, we can conclude that either $V' \cup V$ or $V'' \cup V$ meets $E_l$ in at most a point. Let's assume that this is true of $D^1=V' \cup V$. The minimality of $|D_{l-1} \cap D_{l+1}|$ insures that $D^1$ is in fact a compressing disk for $F^k_p$. We may also isotope $D^1$ so that it misses $D_{l-1}$ entirely, and so that $|D^1 \cap D_{l+1}|<|D_{l-1} \cap D_{l+1}|$. Hence, we have the following chain in $\Gamma (F^k_p)$:
\[D_{l-1}-E_l-D^1-E_l-D_{l+1}\]

We now repeat the construction, using $D^1$ and $D_{l+1}$ to define a disk, $D^2$, such that $D^1 \cap D^2=\emptyset$, and $|D^2 \cap D_{l+1}|<|D^1 \cap D_{l+1}|$, yielding the chain:
\[D_{l-1}-E_l-D^1-E_l-D^2-E_l-D_{l+1}\]

Eventually, this process must terminate with the desired chain. 
\end{proof}

We now proceed to show that our original choice of ${\bf F}$ was not flattened. We do this in two stages. First, for all odd $l$ between 1 and $w-1$, define $B^l$ to be the GHS obtained from $F^k$ by the weak reduction or destabilization, $D_l-E_{l+1}$. For all even $l$, define $B^l$ to be the GHS obtained from $F^k$ by the weak reduction or destabilization, $E_l-D_{l+1}$. Note that since $D_1=D$ and $E_2=E$, it follows that $B^1=F^{k-1}$. Furthermore, Claim \ref{c:isotopy} impies that $B^{w-1}=F^{k+1}$.

Let $A^l$ be a copy of the GHS, $F^k$, for all $l$ between 1 and $w-2$. Now, insert into ${\bf F}$ the sequence of GHSs, $\{B^1,A^1,B^2, A^2, ..., A^{w-2}, B^{w-1}\}$, where $\{F^{k-1},F^k,F^{k+1}\}$ occurred. The new SOG thus defined looks like
\[ ..., F^{k-2}, B^1, A^1, B^2, A^2, ..., B^{w-2}, A^{w-2}, B^{w-1}, F^{k+2}, ...\]

The construction of ${\bf F}'$ from ${\bf F}$ is illustrated in Figure \ref{f:stage1}.

        \begin{figure}[htbp]
        \psfrag{1}{$B^1$}
        \psfrag{2}{$A^1$}
        \psfrag{3}{$B^2$}
        \psfrag{4}{$A^2$}
        \psfrag{8}{$A^{w-2}$}
        \psfrag{9}{$B^{w-1}$}
        \psfrag{a}{$F^{k-2}$}
        \psfrag{b}{$F^{k+2}$}
        \psfrag{f}{$F^k$}
        \psfrag{g}{$F^{k-1}$}
        \psfrag{h}{$F^{k+1}$}
        \vspace{0 in}
        \begin{center}
        \epsfxsize=4 in
        \epsfbox{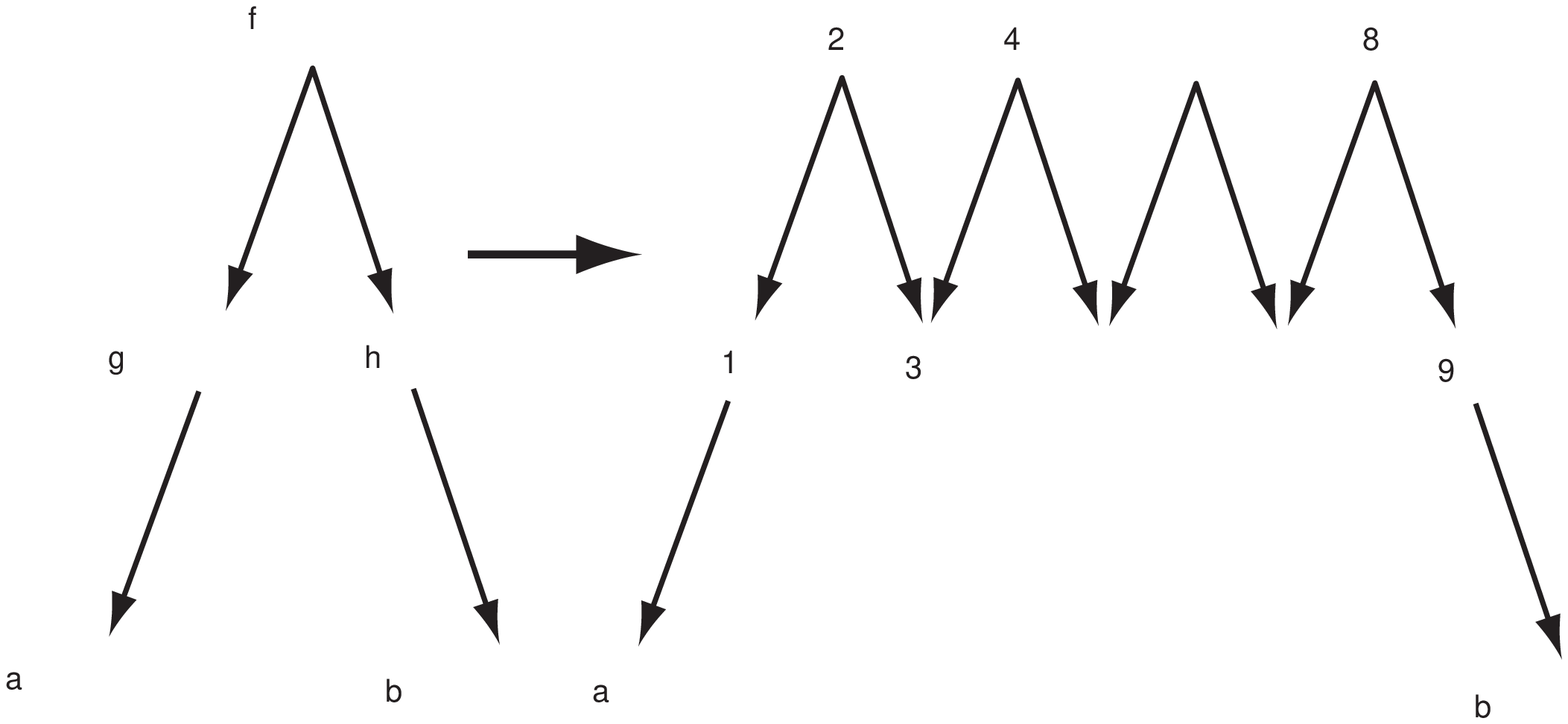}
        \caption{Constructing ${\bf F}'$ from ${\bf F}$.}
        \label{f:stage1}
        \end{center}
        \end{figure}

To see that this defines a SOG, ${\bf F}'$, note that for each odd $l$ (say), $B^l$ is obtained from both $A^{l-1}$ and $A^l$ by the weak reduction or destabilization, $D_l-E_{l+1}$. Note also that each of the $A^l$s that we have inserted becomes a maximal GHS of ${\bf F}'$. 

Now, for each $l$ we perform some operation on ${\bf F}'$ that will result in a new SOG, ``flatter" than ${\bf F}$. We assume that $l$ is even, so that $C_{l \pm 1}=D_{l \pm 1}$, and $C_l=E_l$ (the cases when $l$ is odd are completely symmetric). By Claim \ref{c:disjoint}, we may assume $D_{l-1} \cap D_{l+1}=\emptyset$. Up to symmetry, there are now three cases to consider:

{\bf Case 1.} $D_{l-1} \cap E_l = \emptyset$, and $E_l \cap D_{l+1} = \emptyset$. 

In this case, $D_{l-1}-E_l$ and  $E_l - D_{l+1}$ persist as weak reductions for $B^l$ and $B^{l-1}$, respectively. Performing either of these weak reductions on the corresponding GHS yields the same GHS, which we call $\bar A^{l-1}$. We now replace $A^{l-1}$ with $\bar A^{l-1}$ in ${\bf F}'$. This is illustrated schematically in Figure \ref{f:case1}.

        \begin{figure}[htbp]
        \psfrag{A}{$A^{l-1}$}
        \psfrag{a}{$\bar A^{1-1}$}
        \psfrag{B}{$B^{l-1}$}
        \psfrag{b}{$B^l$}
        \psfrag{1}{$D_{l-1}-E_l$}
        \psfrag{2}{$E_l - D_{l+1}$}
        \vspace{0 in}
        \begin{center}
        \epsfxsize=5 in
        \epsfbox{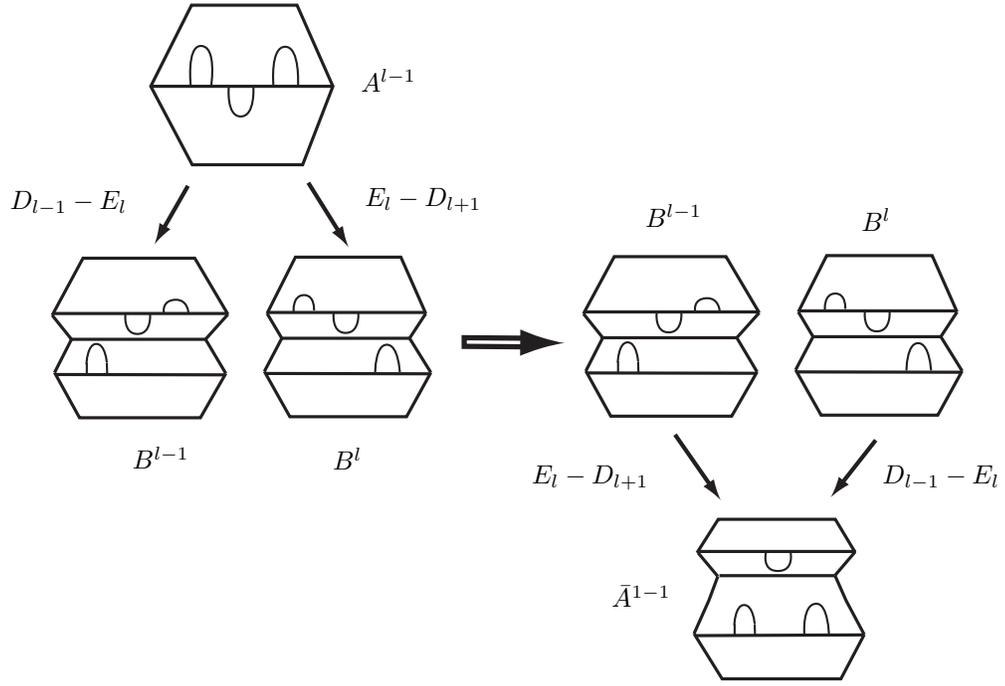}
        \caption{Replacing $A^{l-1}$ with $\bar A^{l-1}$.}
        \label{f:case1}
        \end{center}
        \end{figure}

{\bf Case 2.} $D_{l-1} \cap E_l = \emptyset$, and $|E_l\cap D_{l+1}| = 1$. 

In this case, $E_l - D_{l+1}$ persists as a destabilization for $B^{l-1}$. Performing this destabilization on $B^{l-1}$ yields the GHS, $B^l$. Hence, we can remove $A^{l-1}$ from ${\bf F}'$, and still be left with a SOG. This is illustrated in Figure \ref{f:case2}.

        \begin{figure}[htbp]
        \psfrag{A}{$A^{l-1}$}
        \psfrag{a}{$\bar A^{1-1}$}
        \psfrag{B}{$B^{l-1}$}
        \psfrag{b}{$B^l$}
        \psfrag{1}{$D_{l-1}-E_l$}
        \psfrag{2}{$E_l - D_{l+1}$}
        \vspace{0 in}
        \begin{center}
        \epsfxsize=5 in
        \epsfbox{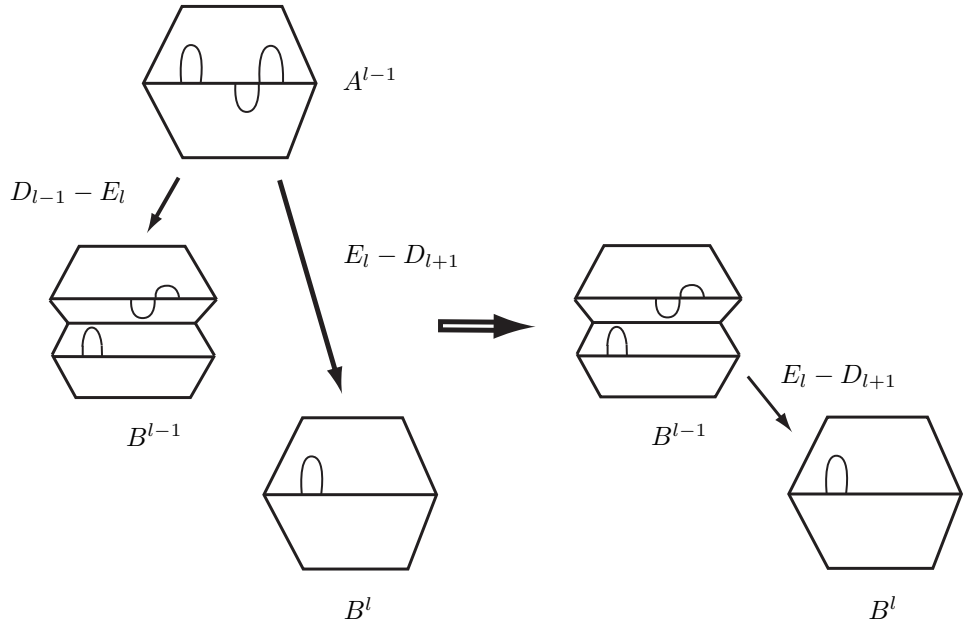}
        \caption{Removing $A^{l-1}$.}
        \label{f:case2}
        \end{center}
        \end{figure}

{\bf Case 3.} $|D_{l-1} \cap E_l| = |E_l\cap D_{l+1}| = 1$. 

Consider the neighborhood of $D_{l-1} \cup E_l \cup D_{l+1}$. This is a ball, which intersects $F^k_p$ in a standardly embedded, twice punctured torus. Hence, the result of the destabilization, $D_{l-1}-E_l$, is isotopic to the result of the destabilization, $E_l-D_{l+1}$. In other words, $B^{l-1}=B^l$, up to isotopy. Hence, we can simply remove $\{B^l,A^l\}$ from ${\bf F}'$, and still be left with a SOG. 

After performing one of the above three operations for each $l$, we eliminate all occurrences of $A^l$. The overall effect is a new SOG in which there are many more maximal GHSs than there were in ${\bf F}$, but where each one is smaller than $F^k$. This shows that $\bf F$ was not flattened, a contradiction.
\end{proof}

We pause here for a moment to note that the proof of Lemma \ref{l:critmax} never used our assumption that $M$ is small. Hence, as an immediate Corollary we obtain the following:

\begin{cor}
If a 3-manifold contains non-isotopic Heegaard splittings of some genus, then it contains a critical surface. 
\end{cor}

This isn't quite the converse of Theorem \ref{t:nocrit}, since we do not conclude that there is a critical {\it Heegaard} surface. For this, we will need our assumption that $M$ is small.

\begin{lem}
\label{l:critsfce}
If $F^k$ is maximal in ${\bf F}$, then $F^k=\{F^k_0, F^k_1, F^k_2\}$, for some critical Heegaard surface, $F^k_1$,
\end{lem}

{\bf Note:} Much of this proof is similar to one which appears in \cite{st:94}. 

\begin{proof}
Suppose that $F^k$ is maximal in ${\bf F}$. By Lemma \ref{l:critmax} there is some odd number, $p$, such that $F^k_p$ is critical in $M^k_p$. We now claim that for each even $i$, $F^k_i$ is incompressible in $M$. Suppose this is not true for $F^k_i$. Then there is a compressing disk, $D$, for $F^k_i$ in $M$. Choose $D$ so as to minimize $w=|D \cap \bigcup \limits _{j\ {\rm even}} F^k_j|$. Let $\alpha$ be an innermost disk of intersection of $D \cap \bigcup \limits _{j\ {\rm even}} F^k_j$, and suppose $\alpha$ lies in $F^k_j$. Then $\alpha$ bounds a subdisk, $D'$, of $D$. If $\alpha$ were inessential on $F^k_j$, then we could swap a subdisk of $F^k_j$ with a subdisk of $D$ to lower $w$, a contradiction. So it must be the case that $D'$ is a compressing disk for $F^k_j$. 

The disk $D'$ lies in either $M^k_{j-1}$ or $M^k_{j+1}$. Without loss of generality, assume the former. By definition, $F^k_{j-1}$ is a Heegaard splitting for $M^k_{j-1}$. However, if $j-1=p$, then $F^k_{j-1}$ is a critical Heegaard splitting for $M^k_{j-1}$, and Corollary \ref{c:essential_boundary} implies that $\partial M^k_{j-1}$ is incompressible in $M^k_{j-1}$. On the other hand, if $j-1 \ne p$, then Lemma \ref{l:oddstrirr} implies $F^k_{j-1}$ is a strongly irreducible Heegaard splitting for $M^k_{j-1}$, and we know that $\partial M^k_{j-1}$ is incompressible in $M^k_{j-1}$ by \cite{cg:87}. 

Since $M$ is assumed to be small, we conclude that there are no non-empty, non-boundary parallel surfaces of even subscript in $F^k$, other than possibly the first or last surface of $F^k$. 

Let $F^k_0$ denote the first surface of $F^k$ and $F^k_m$ the last. Suppose then that $F^k_i$ is boundary parallel, where $i$ is some even number not equal to $0$ or $m$. Let $F_0$ denote the components of $F^k_i$ that are parallel to $F^k_0$, and $F_m$ denote the components of $F^k_i$ that are parallel to $F^k_m$. Assume $F_0$ is non-empty. If $F_m$ is also non-empty, then choose $x \in F_0$, and $y \in F_m$. We claim that there is no path in $M$ connecting $x$ to $y$. If $c$ is such a path, then let $c'$ denote the closure of a component of $c - F^k_i$ which connects $F_0$ to $F_m$. If $c'$ sits in the submanifold of $M$ between $F^k_0$ and $F^k_i$, then, since $c' \cap F_0 \ne \emptyset$, and $F_0$ is parallel into $F^k_0$, the other endpoint of $c'$ must lie on either $F_0$ or $F^k_0$, a contradiction. Similarly, if $c'$ sits in the submanifold of $M$ between $F^k_i$ and $F^k_m$, then, since $c' \cap F_m \ne \emptyset$, and $F_m$ is parallel into $F^k_m$, the other endpoint of $c'$ must lie on either $F_m$ or $F^k_m$, a contradiction. 

We conclude that $F_m =\emptyset$, i.e. that every component of $F^k_i$ is parallel to some component of $F^k_0$. But then this must also be true for $F^k_j$, for {\it all} even $j$ between $0$ and $i$ (since they are all incompressible, and the only such surfaces in a product are boundary parallel). In particular, $F^k_2$ is parallel into $F^k_0$, and so $F^k_1$ is a strongly irreducible Heegaard splitting of a product. As the only such splittings are trivial by \cite{st:93}, and trivial Heegaard splittings do not occur among the surfaces of odd subindex of a GHS, we conclude that $m=2$, i.e. that there are no non-empty surfaces of even subscript in $F^k$, except for possibly $F^k_0$ and $F^k_2$.
\end{proof}

\begin{lem}
\label{l:minHeegaard}
There is a flattened SOG, $\bf F$, whose minimal GHSs are of the form $\{F^k_0, F^k_1, F^k_2\}$, for some strongly irreducible Heegaard surface, $F^k_1$. 
\end{lem}

\begin{proof}
Let $\bf F$ first denote any flattened SOG, and $F^k$ a GHS which is minimal in $\bf F$. If for some odd $i$  $F^k_i$ is not strongly irreducible in $M^k_i$ then there is some edge in $\Gamma(F^k_i)$. We can now form a new GHS, $F^*$, from $F^k$ by doing the corresponding weak reduction or destabilization, and replace $F^k$ in ${\bf F}$ with the sequence $F^k, F^*, F^k$. This replaces one minimal GHS of ${\bf F}$ with a smaller one, without changing any other maximal or minimal GHS. We may now repeat this process until we arrive at a flattened SOG in which all odd surfaces of all minimal GHSs are strongly irreducible.  

Scharlemann and Thompson show in \cite{st:94} that if $F^k_i$ is strongly irreducible in $M^k_i$, for all odd $i$, then $F^k_i$ is incompressible in $M$ for all even $i$. The proof of the Lemma is finished by noting that since $M$ is small it contains no non-boundary parallel incompressible surfaces. Hence, as in the proof of Lemma \ref{l:critsfce}, $F^k$ can contain no non-empty surfaces with even subscript, aside from the first and last. 
\end{proof}

\begin{lem}
If $F^k=\{F^k_0, F^k_1, F^k_2\}$ is maximal in ${\bf F}$, and $F^l=\{F^l_0, F^l_1,F^l_2\}$ is the next (or previous) minimal GHS in ${\bf F}$, then $F^k_1$ is a stabilization of $F^l_1$, $F^k_0=F^l_0$, and $F^k_2=F^l_2$.
\end{lem}

\begin{proof}
Let $\Theta$ denote all of the possibilities for the sequence of GHSs between $F^k$ and $F^l$. Suppose that $\{X^i\}_{i=1}^{n} \in \Theta$ (so that $X^1=F^k$ and $X^n=F^l$). Recall that for each $i$ between $1$ and $n-1$, $X^{i+1}$ is obtained from $X^i$ by either a weak reduction or a destabilization. Note that $X^n$ must be obtained from $X^{n-1}$ by a destabilization, because a weak reduction can never produce a GHS of the form $\{F^l_0, F^l_1,F^l_2\}$ (because the number of surfaces of odd subindex in a GHS which results from a weak reduction is always greater than 1).

We first claim that $F^k_1$ is a stabilization of $F^l_1$. For each $X' \in \Theta$, let $DS(X')$ denote the number of destabilizations which occur after the last weak reduction in $X'$ (If there are no weak reductions in some $X'$, then our claim is proved). By the last sentence of the above paragraph, we know that $DS(X') \ge 1$, for all $X' \in \Theta$. Let $WR(X')$ denote the total number of weak reductions in $X'$. Finally, let $X$ denote the element of $\Theta$ for which the pair $(WR(X),DS(X))$ is minimal, where such pairs are compared lexicographically.

Suppose $X^i$ is the last GHS in $X$ obtained from $X^{i-1}$ by weak reduction. This weak reduction involves compressing disks, $D$ and $E$, for some odd surface, $X^{i-1}_p$, of $X^{i-1}$ such that $D \cap E=\emptyset$. By assumption, $X^{i+1}$ is obtained from $X^i$ by a destabilization. This destabilization involves compressing disks, $D'$ and $E'$, for some odd surface of $X^i$ such that $|D' \cap E'|=1$. Now it is a matter of enumerating all possible cases, and checking the definitions of weak reduction and destabilization. In each case we find that we can either reduce the total number of weak reductions by one, or we can switch the order of the last weak reduction and the next destabilization. Either of these will reduce our complexity, providing a contradiction. We will do the most difficult (and illustrative) case here, and leave the others as an exercise for the reader.

Let $X^{i-1}_D$, $X^{i-1}_E$, and $X^{i-1}_{DE}$ be the surfaces obtained from $X^{i-1}_p$ by compression along $D$, $E$ and both $D$ and $E$. The first case of the definition of a weak reduction is when $X^{i-1}_D \ne X^{i-1}_{p-1}$, and $X^{i-1}_E \ne X^{i-1}_{p+1}$. In this case, $X^i$ is obtained from $X^{i-1}$ by removing the surface, $X^{i-1}_p$, inserting $\{X^{i-1}_D, X^{i-1}_{DE}, X^{i-1}_E\}$ in its place, and reindexing. The matter of reindexing is just for notational convenience, so we may choose to hold off on this for a moment. Hence, for now the surfaces of $X^i$ are $\{..., X^{i-1}_{p-1}, X^{i-1}_D, X^{i-1}_{DE}, X^{i-1}_E, X^{i-1}_{p+1}, ...\}$.

If $D'$ and $E'$ are compressing disks for any surface, $X^{i-1}_j$, where $j$ is an odd number not equal to $p$, then clearly we could have done the destabilization, $D'-E'$ before the weak reduction, $D-E$. Such a switch reduces $DS(X')$, keeping $WR(X')$ fixed, contradicting our minimality assumption. We conclude then that $D'$ and $E'$ are compressing disks for either $X^{i-1}_D$ or $X^{i-1}_E$. Assume the former. 

Now we check the defintion of what it means to peform the destabilizaton, $D'-E'$. The relevant surfaces of $X^i$ are $\{..., X^{i-1}_{p-1}, X^{i-1}_D, X^{i-1}_{DE}, ...\}$. By defintion, $D'$ lies in the submanifold of $M$ cobounded by $X^{i-1}_{p-1}$ and $X^{i-1}_D$, while $E'$ lies in the submanifold cobounded by $X^{i-1}_D$ and $X^{i-1}_{DE}$. But $X^{i-1}_{DE}$ was obtained from $X^{i-1}_D$ by compression along $E$. Hence, it must be the case that $E'=E$. 

To form the GHS, $X^{i+1}$, we are now instructed to create the surfaces, $X^i_{D'}$ and $X^i_{E'}$, which are obtained from $X^{i-1}_D$ by compression along $D'$ and $E'$. Since $E'=E$, it must be that $X^i_{E'}=X^{i-1}_{DE}$. Assume that $X^i_{D'} \ne X^{i-1}_{p-1}$ (if this is not true, then there is another case to check). To form $X^{i+1}$, we are instructed to simply remove $\{X^{i-1}_D, X^{i-1}_{DE}\}$ from $X^i$, and reindex. The resulting GHS now looks like $\{..., X^{i-1}_{p-1}, X^{i-1}_E, X^{i-1}_{p+1}, ...\}$. We now claim that there was a destabilization of $X^{i-1}$ which would have resulted in precisely this GHS. 

Since the surface, $X^i_{D'}$, was obtained from $X^{i-1}_D$ by compression along the red disk, $D'$, and $X^{i-1}_D$ was obtained from $X^{i-1}_p$ by compression along the red disk, $D$, we can identify $D'$ with a red compressing disk (which we will continue to call $D'$) for $X^{i-1}_p$. That is, there is at least one compressing disk for $X^{i-1}_p$ which will become isotopic to $D'$ after we compress along $D$. Choose one, and call it $D'$ also. As $E'=E$, and $|D' \cap E'|=1$, we have $|D' \cap E|=1$. Hence, we can perform the destabilization, $D'-E$, on the GHS, $X^{i-1}$. Let's do this, and see what we get. 

One more time, the first step in performing a destabilization is to form the surfaces, $X^{i-1}_{D'}$, and $X^{i-1}_E$, obtained from $X^{i-1}_p$ by compression along $D'$ and $E$. We have already assumed that $X^{i-1}_E \ne X^{i-1}_{p+1}$. Since neither $X^{i-1}_D$ nor $X^i_{D'}$ were equal to $X^{i-1}_{p-1}$, it must be that $X^{i-1}_{D'} \ne X^{i-1}_{p-1}$. Now, checking the definition of destabilization one last time, to perform $D'-E$ on $X^{i-1}$, we are instructed to remove $X^{i-1}_p$, and in its place insert $X^{i-1}_E$. The resulting GHS now looks like $\{..., X^{i-1}_{p-1}, X^{i-1}_E, X^{i-1}_{p+1}, ...\}$, which we have already seen is precisely the GHS, $X^{i+1}$. What we have shown is that we can reduce the total number of weak reductions by one, again contradicting our minimality assumption on $X$. 

The above argument shows that $F^k_1$ is a stabilization of $F^l_1$. The proof is now complete by noting that stabilization does not change the induced partition of the boundary components of $M$, and so it must be that $F^k_0=F^l_0$, and $F^k_2=F^l_2$.
\end{proof}

By the previous Lemmas, we may assume that for all $k$, $F^k=\{\partial _1 M, F^k_1, \partial _2 M\}$, for some Heegaard surface, $F^k_1$. If $F^k$ is minimal in ${\bf F}$, then $F^k_1$ is strongly irreducible. If $F^k$ is maximal in ${\bf F}$, then $F^k_1$ is critical. Furthermore, for all $k$ either $F^k$ or $F^{k+1}$ is obtained from the other by a destabilization. 

\begin{lem}
There is exactly one maximal GHS in ${\bf F}$. 
\end{lem}

\begin{proof}
Suppose not. Let $F^a$ and $F^c$ be consequetive maximal GHSs of ${\bf F}$, and $F^b$ be the minimal GHS such that $a<b<c$. We assume, without loss of generality, that $c-b \ge b-a$. Note that $F^b$ is obtained from both $F^a$ and $F^{2b-a}$ by a sequence of $b-a$ destabilizations. In other words, the Heegaard splittings, $F^a_1$ and $F^{2b-a}_1$ can be obtained from the Heegaard splitting, $F^b_1$, by doing $b-a$ stabilizations. However, stabilization is unique, so $F^a=F^{2b-a}$. Hence, if we remove the GHSs between $F^{a+1}$ and $F^{2b-a}$ (inclusive) from ${\bf F}$, we are left with a new SOG, with one fewer maximal GHS ($F^a$ will no longer be maximal, unless $c-b=b-a$. In this case, we will have removed $F^b$), where every other maximal GHS has remained unchanged. This contradicts our minimality assumption on ${\bf F}$.
\end{proof}

We summarize all of our results as follows: there is a SOG, ${\bf F}$, of $M$ such that
\begin{enumerate}
    \item $F^1=\{\partial _1 M, F, \partial _2 M\}$,
    \item $F^n=\{\partial _1 M, F', \partial _2 M\}$,
    \item For all $k$, $F^k=\{\partial _1 M, F^k_1, \partial _2 M\}$, for some Heegaard surface, $F^k_1$. If $F^k$ is maximal in ${\bf F}$, then $F^k_1$ is critical.
    \item There is a unique maximal GHS, $F^k$, for ${\bf F}'$. 
\end{enumerate}

We can say all this in a much simpler way: $F^k_1$ is a critical Heegaard splitting of $M$, which is a common stabilization of $F$ and $F'$. Our assumption that ${\bf F}$ is flattened insures that $F^k_1$ is a {\it minimal genus} common stabilization, and hence the proof of Theorem \ref{t:common_stab} is complete.

Now, if we combine Theorem \ref{t:common_stab} with Theorem \ref{t:nocrit}, we obtain the following:

\begin{cor}
A small 3-manifold contains a critical Heegaard surface if and only if it contains non-isotopic Heegaard splittings of some genus. 
\end{cor}

\section{Open Questions}
\label{s:questions}

We conclude with some open questions about critical surfaces. 

\subsection{A metric on the space of strongly irreducible Heegaard splittings}

We now show how our results lead to a natural metric on the space of strongly irreducible Heegaard splittings of a non-Haken 3-manifold. The author believes that it would be of interest to understand this space better. 

First, given a critical surface, $F$, we can define a larger 1-complex, $\Lambda(F)$, that contains $\Gamma(F)$ as follows: the vertices of $\Lambda(F)$ are equivalence classes of loops on $F$, where two loops are considered equivalent if there is an isotopy of $M$ taking $F$ to $F$, and one loop to the other. There is an edge connecting two vertices if there are representatives of the corresponding equivalence classes which intersect in at most a point. Recall that a vertex of $\Gamma(F)$ corresponds to an equivalence class of compressing disks for $F$. Thus, we can identify each vertex of $\Gamma(F)$ with the vertex of $\Lambda(F)$ which corresponds to the boundary of any representative disk. 

Now, suppose $e_1$ and $e_2$ are two edges in $\Gamma (F)$. Define $d(e_1,e_2)$ to be the minimal length of any chain connecting $e_1$ to $e_2$ in $\Lambda(F)$. Now, given two components, $C_1$ and $C_2$, of $\Gamma (F)$, we can define $d(C_1,C_2)$, the {\it distance} between $C_1$ and $C_2$, as $\min \{d(e_1,e_2)|e_1$ is an edge in $C_1$, and $e_2$ is an edge of $C_2\}$.

Finally, suppose $H$ and $H'$ are strongly irreducible Heegaard splittings of a 3-manifold, $M$, and $F$ is their minimal genus common stabilization. If $D-E$ is a destabilization which leads from $F$ to $H$, and $D'-E'$ is a destabilization leading to $H'$, then the results of this paper show that $D-E$ and $D'-E'$ are in different components, $C$ and $C'$, of $\Gamma (F)$. We can therefore define the {\it distance} between $H$ and $H'$ as $d(C,C')$. 

\begin{quest}
Can the distance between strongly irreducible Heegaard splittings be arbitrarily high? If not, is there a bound in terms of the genera of the splittings, or perhaps a universal bound?
\end{quest}

\begin{quest}
Is there an algorithm to compute the distance between two given strongly irreducible Heegaard splittings?
\end{quest}

\begin{quest}
Is there a relationship between the distance between two strongly irreducible Heegaard splittings, and the number of times one needs to stabilize the higher genus one to obtain a stabilization of the lower genus one?
\end{quest}

\begin{quest}
Is there a relationship between the distance between two strongly irreducible Heegaard splittings, and the distances of each individual splitting, in the sense of Hempel \cite{hempel:01}?
\end{quest}

\subsection{The Algorithms}

In \cite{crit2} we produce an algorithm which enumerates isotopy classes of critical surfaces. This ties critical surfaces to several open algorithmic questions in 3-manifold topology. 

\bigskip

\noindent \underline {\it An algorithm to determine stabilization bounds.}  

If one could give an algorithm which takes a critical surface, and returns all of the strongly irreducible Heegaard splittings that it destabilizes to, one may be able to produce an algorithm which takes two Heegaard splittings, and says {\it exactly} how many times they have to be stabilized before they become equivalent. 

\bigskip

\noindent \underline {\it An algorithm to find examples where many stabilizations are required.} 

Since the results of \cite{crit2} give us a way to enumerate critical surfaces, we can look for examples of manifolds with multiple minimal genus Heegaard splittings, but no critical surfaces of genus one higher. This would imply that the minimal genus Heegaard splittings require more than one stabilization to be equivalent. There is currently no known example of this.

\subsection{Topology}

As with any sort of surface in a 3-manifold, it would be nice to have restrictions on how critical surfaces may intersect other submanifolds. Here are some examples of these types of questions:

\bigskip

\noindent \underline {\it The stabilization conjecture.} 

The stabilization conjecture asserts that if we start with two non-isotopic Heegaard splittings and stabilize the higher genus one once, we obtain a stabilization of the other. For non-Haken 3-manifolds, Theorem \ref{t:common_stab} allows us to rephrase this as follows: Suppose $F$ is a critical Heegaard surface, and $D_0-E_0$ and $D_1-E_1$ are edges in different components of $\Gamma (F)$. Let $B_i$ be a ball which contains $D_i$ and $E_i$, such that the genus of $F \cap B_i$ is maximal. If the genus of $F \cap B_0$ is greater than or equal to 2, then the genus of $F \cap B_1$ is 1. 

\bigskip

\noindent \underline {\it Local Detection of critical surfaces.} 

In \cite{scharlemann:97}, Scharlemann proves that under the right hypotheses, the intersection of a strongly irreducible Heegaard splitting with a ball must be a ``nicely embedded" punctured sphere. Is there a restriction on the topology of the intersection of a critical surface with a ball?

\bigskip

\noindent \underline {\it Critical surfaces as index 2 minimal surfaces.} 

In \cite{fhs:83}, Freedman, Hass and Scott show that any incompressible surface can be isotoped to be a least area surface. Such surfaces are minimal surfaces of index 0. In \cite{pr:87}, Pits and Rubinstein show that strongly irreducible surfaces can be isotoped to minimal surfaces of index 1. This motivates us to make the following conjecture:

\begin{cnj}
Any critical surface can be isotoped to be a minimal surface of index 2. 
\end{cnj}

In \cite{crit2}, we prove a Piecewise-Linear analogue of this.

\bibliographystyle{alpha}
\bibliography{critical}

\end{document}